\documentclass[11pt]{article}
\usepackage{graphicx}
\usepackage{color}
\usepackage{amsmath,amssymb,amsthm,amsfonts}
\usepackage{amssymb}
\usepackage{cite}

\usepackage{subfig}
\usepackage{amsmath}
\allowdisplaybreaks[4]
\newtheorem{thm}{Theorem}[section]

\newtheorem{lem}[thm]{Lemma}
\newtheorem{prop}[thm]{Proposition}

\theoremstyle{definition}
\newtheorem{defn}{Definition}[section]
\theoremstyle{remark}

\numberwithin{equation}{section}

\DeclareMathSymbol{\C}{\mathalpha}{AMSb}{"43}

\newcommand{\lam}{\lambda}
\newcommand{\alp}{\alpha}

\newcommand{\R}{{\mathbb{R}}}

\newcommand{\bsub}{\begin{subequations}}
\newcommand{\esub}{\end{subequations}$\!$}

\begin{document}
\title{Refined Limiting Profiles of the Principal Eigenvalue Problems with Large Advection}

\author{ Yujin Guo\thanks{School of Mathematics and Statistics, and Hubei Key Laboratory of
Mathematical Sciences, Central China Normal University, P.O. Box 71010, Wuhan 430079, P. R. China. Email: \texttt{yguo@ccnu.edu.cn}.
}
,\, 	Yuan Lou\thanks{School of Mathematical Sciences, CMA-Shanghai and MOE-LSC, Shanghai Jiao Tong University, Shanghai 200240,
 P. R.	China. Email: \texttt{yuanlou@sjtu.edu.cn}.
},\,
and\, Hongfei Zhang\thanks{School of Mathematics and Statistics, Central China Normal University, P.O. Box 71010, Wuhan 430079,
P. R. China.  Email: \texttt{hfzhang@mails.ccnu.edu.cn}.}
}
\date{\today}
\smallbreak\maketitle
\begin{abstract}
In this paper, we are concerned with the following eigenvalue problem with an advection term:
\begin{equation}\label{0.1}
\left\{
\begin{split}
-\epsilon\Delta \phi-2\alpha\nabla m(x)\cdot\nabla \phi+V(x)\phi&=\lambda \phi\ \ \text{in}\ \ \Omega,\\
\phi&=0\ \ \hbox{on}\ \ \partial\Omega,
\end{split}
\right.
\end{equation}
where $\Omega\subset\mathbb{R}^N~(N\geq1)$ satisfying $\partial\Omega\in C^{2}$ is a bounded domain and contains the origin as an interior point, the constants $\epsilon>0$ and $\alpha>0$ are the diffusive and advection coefficients, respectively, and $m(x)\in C^{2}(\bar{\Omega})$, $V (x)\in C^{\gamma}(\bar{\Omega})~(0<\gamma<1)$ are given functions. We analyze the refined limiting profiles of the principal eigenpair $(\lambda, \phi)$ for (\ref{0.1}) as $\alpha\rightarrow\infty$, which display the visible effect of the large advection on $(\lambda, \phi)$. It expects that our argument is applicable to investigating the refined expansions of the general principal eigenvalue problems.
\end{abstract}

\vskip 0.05truein

\vskip 0.1truein

\vskip 0.1truein
	
\tableofcontents

\section{Introduction}
There have been extensive investigations on the reaction-diffusion models since the last several decades (cf. \cite{ALL,BC,CC2}), which play important roles in studying biological invasions, climate changes and population dynamics.
The principal eigenvalue is a basic concept in the field of reaction-diffusion equations, and it usually provides the important information for the analysis of nonlinear reaction-diffusion equations. In this paper, we consider the following principal eigenvalue problem with an advection term:
\begin{equation}\label{1.1}
\left\{
\begin{split}
-\epsilon\Delta \phi-2\alpha\nabla m(x)\cdot\nabla \phi+V(x)\phi&=\lambda (\alp)\phi\ \  \text{in}\ \ \Omega,\\
\phi &=0\ \  \ \, \hbox{on}\ \ \partial\Omega,
\end{split}
\right.
\end{equation}
where $\Omega$ satisfying $\partial\Omega\in C^{2}$ is a bounded domain in $\mathbb{R}^{N}~(N\geq1)$, and the origin is an interior point of $\Omega$; the constants
$\epsilon>0$ and $\alpha>0$ stand (cf. \cite{PZZ}) for the
diffusive and advection ($i.e.$, drift) coefficients, respectively. Here
$0\leq V (x)\in C^{\gamma}(\bar{\Omega})~(0<\gamma<1)$ and
$m(x)=g(x)|x|^p$, where $p\in\{2\}\cup(3,\infty)$, and $g(x)$ satisfies that
\begin{itemize}
\item[\rm($M$).] $0<c_0:=g(0)\leq g(x)\in C^{2,\gamma_{1}}(\bar{\Omega})~(0<\gamma_{1}<1)$ holds for some positive constant $c_0$,
    $$\Delta g\geq0\ \ \text{and}\ \ x\cdot\nabla g\geq0\ \ \text{in}\ \ \bar{\Omega},$$
    and there exists a small ball $B_{r_0}(0)\subseteq \Omega$ such that
    $$|\nabla g(x)|\leq C|x|^{l},\ \ |\Delta g(x)|\leq C|x|^{l}\ \ \text{in}\ \ B_{r_0}(0),$$
    where $C>0$ and $l>0$ are independent of $r_0>0$.
\end{itemize}
Under above assumptions on $m(x)$ and $V(x)$, similar to  \cite{PZ,PZZ1}, one can get that for any given $\epsilon>0$, (\ref{1.1}) admits a principal eigenvalue $\lambda(\alp)$, which is unique and corresponds to a unique positive eigenfunction $\phi_\alp$, $i.e.$, the so-called unique principal eigenfunction. For convenience, the above mentioned pair $(\lambda(\alp), \phi_\alp)$ is hence called the unique principal eigenpair of (\ref{1.1}).

The principal eigenvalue problem (\ref{1.1}) also connects (cf. \cite{BC,BHN,PZ}) well with the nonlinear propagation phenomena of reaction-diffusion equations. Therefore, it has attracted a lot of analytical investigations on  (\ref{1.1}), see \cite{FS,BNV,BHN,DEF,DF,Friedman,PZZ,Ventcel} and the references therein. Actually, the parabolic versions of (1.1) were   investigated in \cite{LLPZ,LLPZ1}. In addition, the principal eigenvalue problem (\ref{1.1}) with the other type of boundary conditions was analyzed recently in \cite{BHN,CL,CL1,PZ,PZZ,PZZ1} and the references therein. On the other hand, there exist  (eg. \cite{CC,HA,Hut,M,PP}) some interesting works on the asymptotic expansions of the principal eigenpairs for other elliptic problems. Particularly, when the advection term $\nabla m(x)$ is divergence free, $i.e.$, $\Delta m(x)\equiv 0$ in $\Omega$, the authors in \cite{BHN} studied successfully the limiting behavior of $\lam(\alpha)$ for (\ref{1.1}) as $\alpha\rightarrow\infty$.
Motivated by these facts, in this paper we analyze mainly the principal eigenvalue problem (\ref{1.1}) under above general assumptions on $m(x)$. For any given $\epsilon>0$, the main purpose of the present paper is to study the refined limiting profiles of the principal eigenpair for (\ref{1.1}) in the  large advection limit ($i.e.,$ as $\alpha\rightarrow\infty$), and the case of small diffusions ($i.e.,$ $\epsilon\rightarrow0$) is however left to the companion work \cite{Guo}.

We remark that the limiting behavior of the principal eigenvalue problems was studied recently in  \cite{BHN,DEF,DF,Friedman} by the first integral and some other methods, which however cannot be applicable to the principal eigenvalue problem  (\ref{1.1}) under general assumptions on $m(x)$. Moreover, the refined limiting profiles of solutions for nonlinear elliptic equations were also analyzed in \cite{Grossi,GLW} and the references therein by employing the non-degeneracy of the associated linearized problems, which method seems restricted unfortunately to the nonlinear elliptic problems. To achieve the main purpose of the present paper, we therefore need to investigate new ideas of analyzing the principal eigenpair $(\lambda(\alpha), \phi_{\alpha})$ for the linear problem (\ref{1.1}) under general assumptions on $m(x)$. Roughly speaking, in this paper we shall derive several different Pohozaev type identities and employ the $L^{2}$--constraint condition of (\ref{1.1}) to analyze rigorously  the refined limiting profiles of the unique principal eigenpair for (\ref{1.1}) as $\alpha\rightarrow\infty$, where $\epsilon>0$ is fixed. It expects that the argument developed in the present paper can be employed to handle with  the refined expansions of the general principal eigenvalue problems.

Following \cite{FS,BHN,PZZ}, we  note that the principal eigenvalue problem $\lam (\alpha)$ of $(\ref{1.1})$ can be equivalently characterized by
\begin{equation}\label{2.1}
\begin{split}
\lambda(\alpha)=&\inf_{\{\phi\in H_{0}^{1}(\Omega),~\phi\not\equiv 0\}}
\frac{\int_{\Omega}e^{\frac{2\alpha}{\epsilon}m(x)}(\epsilon|\nabla \phi|^{2}+V(x)\phi^{2})dx}{\int_{\Omega}e^{\frac{2\alpha}{\epsilon}m(x)}\phi^{2}dx}\\
=&\inf_{\{u\in H_{0}^{1}(\Omega),~\int_{\Omega}u^{2}dx=1\}}
\Big\{\int_{\Omega}\epsilon\big|\nabla u-\frac{\alpha}{\epsilon} u
\nabla m\big|^2dx+\int_{\Omega}V(x)u^2dx\Big\}\\
=&\inf_{\{u\in H_{0}^{1}(\Omega),~\int_{\Omega}u^{2}dx=1\}}\Big\{\int_{\Omega}\epsilon|\nabla u|^2dx+\int_{\Omega}\Big[\frac{\alpha^{2}}{\epsilon}|\nabla m|^{2}
\\
&\qquad\qquad\qquad\qquad\quad+\alpha\Delta m+V(x)\Big]u^2dx\Big\},
\end{split}
\end{equation}
where the second variational characterization is derived by substituting  $\phi=e^{-\frac{\alpha}{\epsilon}m}u$.
Moreover, if $u_{\alpha}$ is a minimizer of $(\ref{2.1})$, then it follows from the variational theory that $u_{\alpha}$ satisfies the elliptic equation
\begin{equation}\label{1.6}
-\epsilon\Delta u_{\alpha}+\Big[\frac{\alpha^{2}}{\epsilon}|\nabla m|^{2}+\alpha\Delta m+V(x)\Big]u_{\alpha}
=\lambda(\alpha)u_{\alpha}
\ \ \text{in}\ \, \Omega.
\end{equation}
By the strong maximum principle, one can also derive from $(\ref{2.1})$ and (\ref{1.6}) that $u_{\alpha}$ does not change sign in $\Omega$, which further implies that the principal eigenfunction of $(\ref{1.1})$ does not change sign in $\Omega$ yet. 

Following the above equivalence between $(\ref{1.1})$ and (\ref{2.1}), suppose $(\lambda(\alpha), \phi_{\alpha})$ is the unique principal eigenpair of (\ref{1.1}), then it means equivalently throughout the present whole paper that $(\lambda(\alpha), u_{\alpha})$ satisfies  $(\ref{2.1})$ and hence (\ref{1.6}). Stimulated by (\ref{1.6}), we now introduce the following constraint variational problem: for any fixed $\epsilon>0$ and $p\in\{2\}\cup(3,\infty)$,
\begin{equation}\label{2.3}
\hat{\lambda}:=\inf_{\{u\in\mathcal{H},~\int_{\mathbb{R}^{N}}u^{2}dx=1\}}
\int_{\mathbb{R}^{N}}\Big[\epsilon|\nabla u|^2+\Big(\frac{p^{2}c_{0}^{2}}{\epsilon}|x|^{2p-2}
+pc_{0}(N+p-2)|x|^{p-2}\Big)u^2\Big]dx,
\end{equation}
where  $c_{0}>0$ is as in ($M$),
and the space $\mathcal{H}$ is defined by
\begin{equation}\label{2.5}
\mathcal{H}:=\Big\{u\in H^{1}(\mathbb{R}^{N}): \int_{\mathbb{R}^{N}}|x|^{2p-2}u^{2}dx<\infty\Big\},
\end{equation}
together with the norm
$$\|u\|_{\mathcal{H}}=\Big\{\int_{\mathbb{R}^{N}}\big[|\nabla u|^{2}+(1+|x|^{2p-2})u^{2}\big]dx\Big\}^{\frac{1}{2}}.$$
We shall prove in Lemma \ref{lem:2.1}  below that $\hat{\lambda}>0$, and $\hat{\lambda}$ admits a unique positive minimizer $\hat{u}(x)$, which is also radially symmetric in $|x|$.

\subsection{Main results}

In this subsection, we introduce the main results of the present paper. Associated to the unique positive minimizer $\hat{u}$ of (\ref{2.3}), for convenience we denote $\hat{\psi}_{i}(x)\in C^{2}(\mathbb{R}^{N})\cap L^\infty(\mathbb{R}^{N})$ to be the unique solution of
\begin{equation}\label{1.121}
\left\{
\begin{split}
 -\epsilon\Delta \hat{\psi}_{i}+\Big[\frac{p^{2}c_{0}^{2}}{\epsilon}|x|^{2p-2}
+pc_{0}(N+p-2)|x|^{p-2}-\hat{\lambda}\Big]\hat{\psi}_{i}
=\hat{F}_{i}(x)\ \ \text{in}\ \  \mathbb{R}^{N},\\
\int_{\mathbb{R}^{N}}\hat{\psi}_{i}\hat{u}dx=0,\ \ i=1,2,\quad\qquad\qquad\qquad\qquad\qquad
\end{split}\right.
\end{equation}
where $\hat{F}_{i}(x)$ satisfies
\begin{equation}\label{1.191}
\hat{F}_{i}(x)=\left\{
\begin{split}-\big(x\cdot\nabla V(0)\big)\hat{u},\ \ \qquad\qquad\qquad &\text{if}\ \ i=1,\\
\sum_{|\sigma|=2}\frac{D^{\sigma}V(0)}{\sigma!}\Big(-x^{\sigma}
+\int_{\mathbb{R}^{N}}x^{\sigma}
\hat{u}^{2}dx\Big)\hat{u},\ \
&\text{if}\ \ i=2.
\end{split}\right.
\end{equation}
Here $\sigma=(\sigma_{1},\cdots,\sigma_{N})$ is a multiple index with nonnegative integers $\sigma_{1},\cdots,\sigma_{N}$, and   $$|\sigma|:=\sigma_{1}+\cdots+\sigma_{N},\ \ \sigma!:=\sigma_{1}!\cdots\sigma_{N}!,\ \  x^{\sigma}:=x_{1}^{\sigma_{1}}\cdots x_{N}^{\sigma_{N}},$$
provided that $ x=(x_{1},\cdots,x_{N})\in \R^N$.
We can get from Lemma \ref{lem:A.1} below the existence and uniqueness of $\hat{\psi}_{i}(x)$ for $i=1,2$, due to the facts that $\int_{\mathbb{R}^{N}}\hat{F}_{i}(x)\hat{u}(x)dx=0$ holds for $i=1,2$, and $\hat{u}>0$ is radially symmetric in $|x|$.

Under general assumptions on $m(x)$ and $V(x)$, in this paper we first establish the following refined limiting profiles of the unique principal eigenpair for (\ref{1.1}) as $\alpha\rightarrow\infty$,
no matter whether the minimum points of $V(x)$ and $m(x)$  coincide with each other or not.


\begin{thm}\label{cor 1.2}
Suppose $0\leq V (x)\in C^{\gamma}(\bar{\Omega})\cap C^{2}(B_{r_{1}}(0))$, where $0<\gamma<1$ and $B_{r_{1}}(0)\subseteq\Omega$ is a small domain for some $r_{1}>0$, and assume $m(x)=g(x)|x|^p\geq0$, where $p\in\{2\}\cup(3,\infty)$ and $g(x)$ satisfies ($M$) for some $l>3$. Then for any fixed $\epsilon>0$, the principal eigenpair $(\lambda(\alpha), u_{\alpha})$ of (\ref{1.1}) satisfies
\begin{equation}\label{3.341}
\lambda(\alpha)=\alpha^{\frac{2}{p}}\hat{\lambda}+
V(0)
+\alpha^{-\frac{2}{p}}\sum_{|\sigma|=2}\frac{D^{\sigma}V(0)}{\sigma!}\int_{\mathbb{R}^{N}}
x^{\sigma}\hat{u}^{2}dx
+o(\alpha^{-\frac{2}{p}})
\ \ \text{as}\ \ \alpha\rightarrow\infty,
\end{equation}
and
\begin{equation}\label{3.361}
\begin{split}
\tilde{w}_{\alpha}(x):=&\alpha^{-\frac{N}{2p}}u_{\alpha}
(\alpha^{-\frac{1}{p}}x)
\\ =&\hat{u}(x)+\alpha^{-\frac{3}{p}}\hat{\psi}_{1}(x)+\alpha^{-\frac{4}{p}}\hat{\psi}_{2}(x)
+o(\alpha^{-\frac{4}{p}})\ \ \text{in}\ \ \mathbb{R}^{N}
\ \ \text{as}\ \ \alpha\rightarrow\infty,
\end{split}
\end{equation}
where $(\hat{\lambda}, \hat{u})$ is the unique principal eigenpair of (\ref{2.3}), $u_\alpha(x)\equiv0$ in $\mathbb{R}^{N}\backslash\Omega$, and $\hat{\psi}_{1}$, $\hat{\psi}_{2}$ are given uniquely by (\ref{1.121}) and (\ref{1.191}).
\end{thm}

There are several comments on Theorem \ref{cor 1.2} in order. Firstly, one can note  that the leading term of the principal
eigenvalue (\ref{3.341}) was essentially derived earlier in \cite{BHN}, which was however restricted to the divergence-free
case of $\nabla m(x)$ for (\ref{1.1}), $i.e.$, $\Delta m(x)\equiv 0$ in $\Omega$. As a completion of \cite{BHN}, Theorem \ref{cor 1.2} is fortunately
applicable to the general case where $\nabla  m(x)$ is not divergence-free, $i.e.$, $\Delta m(x)\not\equiv 0$ in $\Omega$, and more importantly, it also presents the higher
order asymptotic behavior of the principal eigenpair for (\ref{1.1}) as $\alpha\rightarrow\infty$. As far as we know, this
seems the first work on the refined limiting profiles of the principal eigenvalue problems. Secondly, one can observe from
(\ref{3.341}) that the profile of  $V(x)$ at the unique minimum point ($i.e.,$ the origin) of $m(x)$ affects strongly the refined
limiting profiles of the principal  eigenpair $(\lambda(\alpha), u_{\alpha})$ for (\ref{1.1}) as $\alpha\rightarrow\infty$.
Moreover, the proof of Theorem \ref{cor 1.2} shows that  the mass of the principal eigenfunction $u_{\alpha}$ concentrates near
the unique minimum point ($i.e.,$ the origin) of  $m(x)$  as $\alpha\rightarrow\infty$, see Proposition \ref{prop2.3} for more details.
Thirdly, the argument of Theorem
\ref{cor 1.2} can yield theoretically more terms for the asymptotic expansions of the principal eigenpair for (\ref{1.1}) as
$\alpha\rightarrow\infty$, provided that one imposes extra assumptions on $V(x)$ and $m(x)$. For simplicity, we leave this generalization to the interested
reader.

We next follow three steps to explain roughly the general strategy of proving Theorem \ref{cor 1.2}, which can be employed to handle with  the refined expansions of the general principal eigenvalue problems. As the first step, we
first derive the least energy estimates of (\ref{1.1})  as $\alpha\rightarrow\infty$, for which the key point is to explore a suitable trial function. This then helps us establish in Lemma \ref{lem:2.2}   the
following limiting behavior of  the principal eigenvalue $\lambda(\alpha)$ for (\ref{1.1}):
\begin{equation}\label{1.4}
\alpha^{-\frac{2}{p}}\lambda(\alpha)\rightarrow\hat{\lambda}>0\ \ \text{as}\ \, \alpha\rightarrow\infty,
\end{equation}
where $\hat{\lambda}$ is given by (\ref{2.3}). Applying (\ref{1.4}), we shall further prove in  Proposition \ref{prop2.3} that the principal eigenfunction $u_{\alpha}$ of (\ref{1.1})
satisfies
\begin{equation}\label{1.5}
\begin{split}
\hat{w}_\alpha(x):=&\alpha^{-\frac{N}{2p}}u_{\alpha}
(\alpha^{-\frac{1}{p}}x+x_\alpha)\\
\to&\,\hat{u}(x)>0\ \ \text{strongly in}\ \ H^{1}(\mathbb{R}^{N})\cap
L^{\infty}(\mathbb{R}^{N})
\ \ \text{as}\ \ \alpha\rightarrow\infty,
\end{split}
\end{equation}
where $\hat{u}$ is the unique positive minimizer of (\ref{2.3}), and
$x_\alpha\in \R^N$ denotes the unique global maximum point of $u_\alpha$ as $\alpha\rightarrow\infty$
and satisfies $\lim_{\alpha\rightarrow\infty}\alpha^{\frac{1}{p}}|x_\alpha|=0$. One can note from Section 2 that the limit behavior of (\ref{1.4}) and (\ref{1.5}) holds essentially for any $p\ge 2$.

The second step of proving Theorem \ref{cor 1.2} is to prove
Lemmas \ref{lem3.1} and \ref{lem3.21} on the refined asymptotic expansions of $(\lambda(\alpha), u_{\alpha})$ as $\alpha\rightarrow\infty$ by deriving several different Pohozaev type identities and employing the $L^{2}$--constraint condition of (\ref{1.1}).
More precisely, the second-order asymptotic expansion of $\lambda(\alpha)$ as
$\alpha\rightarrow\infty$ is obtained in Lemma \ref{lem3.1} by deriving the key identity (\ref{3.3}).
Although the higher order asymptotic expressions of solutions for nonlinear elliptic equations were also explored, see \cite{Grossi,GLW} and the references therein, by applying the non-degeneracy of the associated linearized problems, this approach cannot be applied unfortunately to the linear elliptic problem  (\ref{1.1}).
To prove Lemma \ref{lem3.21} on the refined asymptotic expansion of $ u_{\alpha}$, inspired by (\ref{1.5}), we hence consider the function
\begin{equation}\label{1.13}
\eta_{\alpha}(x):=\tilde{w}_{\alpha}(x)
-\Big(\int_{\mathbb{R}^{N}}\tilde{w}_{\alpha}\hat{u}dx\Big)\hat{u}(x),\ \    \tilde{w}_{\alpha}(x)=\alpha^{-\frac{N}{2p}}u_{\alpha} (\alpha^{-\frac{1}{p}}x),
\end{equation}
where $\int_{\mathbb{R}^{N}}\tilde{w}_{\alpha}\hat{u}dx\to 1$ as $\alp \to\infty$ in view of (\ref{1.5}).
We emphasize that, instead of the constant $1$, the varying coefficient $\int_{\mathbb{R}^{N}}\tilde{w}_{\alpha}\hat{u}dx$ in (\ref{1.13}) is used to ensure the crucial orthogonality between $\eta_{\alpha}(x)$ and $\hat{u}(x)$. By employing fully the $L^{2}$--constraint condition of (\ref{1.1}), we further derive in Lemma \ref{lem3.2} the refined estimate of $\eta_{\alpha}(x)$ as $\alp \to\infty$, which finally helps us to establish Lemma \ref{lem3.21}.


The third step of proving Theorem \ref{cor 1.2} is to establish (\ref{3.341}) and (\ref{3.361}).  By applying $L^2-$constraint conditions of $u_\alpha$, together with (\ref{3.57}), the key point of this step is to derive the important Pohozaev type identity (\ref{3.56}) in terms of the parameter $\theta_{\alpha}$. The parameter $\theta_{\alpha}$ is then identified in (\ref{3.61}), based on which we finally derive  (\ref{3.341}) and (\ref{3.361}). For the detailed proof of  Theorem \ref{cor 1.2}, we refer the reader to Sections 2 and 3.

One can note that Theorem \ref{cor 1.2} holds true for the case where $V(x)$ is smooth enough in the neighborhood of the
origin. For this reason, we are next concerned with the refined limiting profiles of the principal eigenpair for (\ref{1.1}) as $\alpha\rightarrow\infty$, provided that the weaker regularity of $V(x)$ holds in the neighborhood of the origin.  Towards this purpose, we now introduce the homogeneous functions.

\begin{defn}\label{1.7} The function
$h(x)\geq0$ in $\R^{N}$ is called homogeneous of degree
$q\in \R^{+}$ (about the origin), if there exists some constant $q>0$ such that
\begin{equation}\label{v1}
h(tx)=t^{q}h(x)\ \ \hbox{in}\ \R^{N}\ \ \mbox{for any}\ \ t>0.
\end{equation}
\end{defn}

\noindent
Following \cite[Remark 3.2]{Grossi}, the above definition implies that the homogeneous function
$h(x)\in C(\R^{N})$ of degree $q>0$ satisfies
$$0\leq h(x)\leq C |x|^{q} \ \ \mbox{in}\,\ \R^{N}, \ \ \mbox{where}\,\ C=\max\limits_{x\in \partial B_{1}(0)}h(x)>0.$$
Typically, $h(x)=|x|^q$ is homogeneous of degree $q>0$ about the origin.
Moreover, if the homogenous function $h(x)$ satisfies
$\lim_{|x|\rightarrow\infty}h(x)=+\infty$, then $x=0$ is the unique minimum point of $h(x)$ in $\R^{N}$.

Suppose $(\hat{\lambda}, \hat{u})$ is the unique principal eigenpair of (\ref{2.3}). We next assume that
\begin{equation}\label{1.11A}
\begin{split}
& 0\leq \hat{h}(x)\in C^{\gamma}(\R^{N})~(0<\gamma<1)\ \ \text{is homogeneous of degree}\ \ \hat{q}>0,\\
&\text{and satisfies}\ \lim_{|x|\rightarrow\infty}\hat{h}(x)=\infty.
\end{split}
\end{equation}
We denote $\hat{\psi}_{i}(x)\in C^{2}(\mathbb{R}^{N})\cap L^\infty(\mathbb{R}^{N})$ ($i=3,4$) to be the unique solution of
\begin{equation}\label{1.12}
\left\{
\begin{split}
 -\epsilon\Delta \hat{\psi}_{i}+\Big[\frac{p^{2}c_{0}^{2}}{\epsilon}|x|^{2p-2}
+pc_{0}(N+p-2)|x|^{p-2}-\hat{\lambda}\Big]\hat{\psi}_{i}
=\hat{F}_{i}(x) \ \ \mbox{in} \ \, \mathbb{R}^{N},\\
\int_{\mathbb{R}^{N}}\hat{\psi}_{i}\hat{u}dx=0,\ \ i=3,4,\qquad\qquad\qquad\qquad\qquad
\end{split}\right.
\end{equation}
where $\hat{F}_{i}(x)$ satisfies
\begin{equation}\label{1.19}
\hat{F}_{i}(x)=
\left\{\begin{split}
-\hat{h}(x)\hat{u}+\Big(\int_{\mathbb{R}^{N}}\hat{h}(x)\hat{u}^{2}dx\Big)
\hat{u},\ \ \qquad\qquad\qquad&\text{if}\ \ i=3,\\
-\hat{h}(x)\hat{\psi}_{3}+
\Big(\int_{\mathbb{R}^{N}}\hat{h}(x)\hat{u}\hat{\psi}_{3}dx\Big)\hat{u}
+\Big(\int_{\mathbb{R}^{N}}\hat{h}(x)\hat{u}^{2}dx\Big)\hat{\psi}_{3},\ \
&\text{if}\ \ i=4.
\end{split}\right.
\end{equation}
Similar to (\ref{1.121}) and (\ref{1.191}), we note that the existence and uniqueness of $\hat{\psi}_{i}(x)$ ($i=3,4$) hold in view of Lemma \ref{lem:A.1} below.

Following
the above notations, the second main result of the present paper can be stated as the following theorem.

\begin{thm}\label{thm:1.3}
Suppose that $0\leq V (x)\in C^{\gamma}(\bar{\Omega})~(0<\gamma<1)$ satisfies $V(x)=\hat{h}(x)$ in some small domain $B_{r_{2}}(0)\subseteq\Omega$ for some $r_{2}>0$, where $\hat{h}(x)$ satisfies (\ref{1.11A}) for some $\hat{q}>0$, and assume that $m(x)=g(x)|x|^p\geq0$, where $p\in\{2\}\cup(3,\infty)$ and $g(x)$ satisfies ($M$) for some  $l>2\hat{q}+3$.
Then for any fixed $\epsilon>0$, the principal eigenpair $(\lambda(\alpha), u_{\alpha})$ of (\ref{1.1}) satisfies
\begin{equation}\label{3.34}
\begin{split}
\lambda(\alpha)=&\alpha^{\frac{2}{p}}\hat{\lambda}+
\alpha^{-\frac{\hat{q}}{p}}\int_{\mathbb{R}^{N}}\hat{h}(x)\hat{u}^{2}dx+
\alpha^{-\frac{2\hat{q}+2}{p}}\int_{\mathbb{R}^{N}}\hat{h}(x)\hat{u}\hat{\psi}_{3}dx
\\&+o(\alpha^{-\frac{2\hat{q}+2}{p}})
\ \ \text{as}\ \ \alpha\rightarrow\infty,
\end{split}
\end{equation}
and
\begin{equation}\label{3.36}
\begin{split}
\tilde{w}_{\alpha}(x):=&\alpha^{-\frac{N}{2p}}u_{\alpha}
(\alpha^{-\frac{1}{p}}x)\\
=&\hat{u}(x)+\alpha^{-\frac{\hat{q}+2}{p}}\hat{\psi}_{3}(x)
+\alpha^{-\frac{2\hat{q}+4}{p}}\Big[\hat{\psi}_{4}(x)
-\frac{1}{2}\Big(\int_{\mathbb{R}^{N}}\hat{\psi}_{3}^{2}dx\Big)\hat{u}(x)\Big]
\\
&+o(\alpha^{-\frac{2\hat{q}+4}{p}})\ \ \text{in}\ \ \mathbb{R}^{N}
\ \ \text{as}\ \ \alpha\rightarrow\infty,
\end{split}
\end{equation}
where $(\hat{\lambda}, \hat{u})$ is the unique principal eigenpair of (\ref{2.3}), $u_\alpha(x)\equiv0$ in $\mathbb{R}^{N}\backslash\Omega$, and $\hat{\psi}_{3}$, $\hat{\psi}_{4}$ are given uniquely by (\ref{1.12}) and (\ref{1.19}).
\end{thm}

We shall only sketch the proof of Theorem \ref{thm:1.3} for simplicity, since it follows from the general strategy  of proving Theorem \ref{cor 1.2}. Similar to Theorem \ref{cor 1.2}, if the further information of  $g(x)$ near the origin is imposed for Theorem \ref{thm:1.3}, then one can expect the higher-order  asymptotic expansions of $(\lambda(\alpha), u_{\alpha})$ as $\alpha\rightarrow\infty$. However, this process does not require the additional assumptions on $V(x)$, which is different from that of Theorem \ref{cor 1.2}.

This paper is organized as follows.  For a general class of $V(x)$, in Section 2 we study the limiting behavior of the principal eigenpair $(\lambda(\alpha),u_{\alpha})$ for (\ref{1.1}) as $\alpha\rightarrow\infty$. The asymptotic expansions of $(\lambda(\alpha),u_{\alpha})$ as $\alpha\rightarrow\infty$ are addressed in Section 3, based on which we then complete in Subsection 3.1 the proof of Theorem \ref{cor 1.2}.  The proof of Theorem \ref{thm:1.3} is finally given in Section 4.

\section{Limiting Behavior of  $(\lambda(\alpha), u_{\alpha})$ as $\alpha\to\infty$}
In this section, we shall establish the limiting behavior of the principal eigenpair $(\lambda(\alpha), u_{\alpha})$ for (\ref{1.1}) as $\alpha\to\infty$. To this end, we start with the following analytical properties of the problem $\hat{\lambda}$ defined by (\ref{2.3}), whose proof is given briefly for the reader's convenience.

\begin{lem}\label{lem:2.1}
For any fixed $\epsilon>0$ and $p\geq2$, consider the problem $\hat{\lambda}$ defined by (\ref{2.3}). Then we have
\begin{enumerate}
\item  $\hat{\lambda}>0$ and (\ref{2.3}) admits a unique positive minimizer $\hat{u}$ in $\mathcal{H}$, which is radially symmetric and strictly decreasing in $|x|>0$.

\item There exist large positive constants $R_{0}>0$ and $C>0$ such that
\begin{equation}\label{2.9}
\hat{u}(x)\leq Ce^{-|x|},\ \ |\nabla\hat{u}(x)|\leq Ce^{-\frac{|x|}{2}}\ \ \text{uniformly for}\ \ |x|\geq R_{0}.
\end{equation}
\end{enumerate}
\end{lem}

\noindent\textbf{Proof.}
1. For any fixed $\epsilon>0$ and $p\geq2$, it is standard to check that (\ref{2.3}) admits at least one minimizer $\hat{u}$ in $\mathcal{H}$. Moreover, it yields from \cite[Theorem 11.8]{Lieb1} that $\hat{u}$ can be chosen to be  strictly positive, and more importantly, $\hat{u}>0$ must be unique in $\mathcal{H}$. By the standard rearrangement theory (see \cite{BBJ,Lieb1}), one can prove that the positive minimizer $\hat{u}$ is
radially symmetric and decreasing in $|x|>0$.
Additionally, using the variational theory, we obtain that $\hat{u}$ solves the following Euler-Lagrange equation
\begin{equation}\label{2.10}
-\epsilon\Delta \hat{u}+\Big[\frac{p^{2}c_0^{2}}{\epsilon}|x|^{2p-2}+ pc_0(N+p-2)|x|^{p-2}\Big]\hat{u}=\hat{\lambda}\hat{u}\ \ \text{in}\ \ \mathbb{R}^{N}.
\end{equation}
By the standard elliptic regularity theory, it then yields from (\ref{2.10}) that $\hat{u}\in C^{2}(\mathbb{R}^{N})$. It thus follows from \cite[Theorem 2]{LN} that $\hat{u}$ is strictly decreasing in $|x|>0$. Finally, since (\ref{2.3}) admits a minimizer $\hat{u}>0$, one can get directly from (\ref{2.3}) that $\hat{\lambda}>0$.

2. We derive from (\ref{2.10}) that there exists a sufficiently large positive constant $R_{0}>0$ such that
\begin{equation*}
-\Delta \hat{u}+\hat{u}<0\ \ \text{uniformly for}\ \ |x|\geq R_{0}.
\end{equation*}
By the comparison principle, we then obtain from above that
\begin{equation}\label{2.11M}
\hat{u}(x)\leq Ce^{-|x|}\ \ \text{uniformly for}\ \ |x|\geq R_{0}.
\end{equation}
Moreover, applying the local elliptic estimate (cf. \cite[(3.15)]{GT}), we derive from (\ref{2.10}) and (\ref{2.11M}) that
\begin{equation*}\label{2.44}
|\nabla\hat{u}(x)|\leq Ce^{-\frac{|x|}{2}}\ \ \text{uniformly for}\ \ |x|\geq R_{0}.
\end{equation*}
This completes the proof of Lemma \ref{lem:2.1}.
\qed

We comment that the exponential decay (\ref{2.9}) is not optimal, but it is enough for the analysis of the present whole paper. Applying Lemma \ref{lem:2.1}, we next analyze the following limiting behavior of the principal eigenvalue $\lambda(\alpha)$ for $(\ref{1.1})$ as $\alpha\rightarrow\infty$, where the  potential $V(x)$ is more general than those of Theorems \ref{cor 1.2} and \ref{thm:1.3}.

\begin{lem}\label{lem:2.2}
Suppose $0\leq V (x)\in C^{\gamma}(\bar{\Omega})~(0<\gamma<1)$, and assume  $m(x)=g(x)|x|^p\geq0$, where $p\geq2$ and $g(x)$ satisfies ($M$) with $l>0$.
Then for any fixed $\epsilon>0$, the principal eigenvalue $\lambda(\alpha)$ of $(\ref{1.1})$ satisfies
\begin{equation}\label{1.2}
\lim_{\alpha\rightarrow\infty}\frac{\lambda(\alpha)}{\alpha^{\frac{2}{p}}}=\hat{\lambda}>0,
\end{equation}
where $\hat{\lambda}$ is defined by (\ref{2.3}).
\end{lem}

\noindent\textbf{Proof.}  Under the assumptions of Lemma \ref{lem:2.2}, suppose $(\lambda(\alpha), u_{\alpha})$ is the unique principal eigenpair of (\ref{1.1}) as $\alpha\to\infty$.
We then deduce from ($M$), (\ref{2.1}) and (\ref{2.3}) that
\begin{align}
\lambda(\alpha)=&\int_{\Omega}\epsilon|\nabla u_{\alpha}|^2dx+\int_{\Omega}\Big[\frac{\alpha^{2}}{\epsilon}|\nabla m|^{2}
+\alpha\Delta m+V(x)\Big]u_{\alpha}^2dx\nonumber\\
=&\int_{\Omega}\epsilon|\nabla u_{\alpha}|^2dx+\int_{\Omega}\Big\{\alpha|x|^p\Delta g+\alpha\big[2p(x\cdot\nabla g)+p(N+p-2)g(x)\big]|x|^{p-2}\nonumber\\
&+\frac{\alpha^{2}}{\epsilon}\Big[|x|^{2p}|\nabla g|^{2}+p|x|^{2p-2}(x\cdot\nabla g^{2})+p^{2}g^{2}(x)|x|^{2p-2}\Big]
+V(x)\Big\}u_{\alpha}^2dx\nonumber\\
\geq&\int_{\mathbb{R}^{N}}\epsilon|\nabla u_{\alpha}|^2dx+\int_{\mathbb{R}^{N}}\Big[\frac{\alpha^{2}p^{2}}{\epsilon}g^{2}(x)|x|^{2p-2}
+\alpha p(N+p-2)g(x)|x|^{p-2}\Big]u_{\alpha}^2dx\nonumber\\
=&\alpha^{\frac{2}{p}}\Big\{\int_{\mathbb{R}^{N}}\epsilon|\nabla \hat{\nu}_{\alpha}|^{2}dx+\int_{\mathbb{R}^{N}}\Big[p(N+p-2)g(\alpha^{-\frac{1}{p}}x)|x|^{p-2}\label{2.12}\\
&\quad\ \ +\frac{p^{2}}{\epsilon}g^{2}(\alpha^{-\frac{1}{p}}x)
|x|^{2p-2}
\Big]\hat{\nu}_{\alpha}^2dx\Big\}\nonumber\\
\geq&\alpha^{\frac{2}{p}}\Big\{\int_{\mathbb{R}^{N}}\epsilon|\nabla \hat{\nu}_{\alpha}|^{2}dx+\int_{\mathbb{R}^{N}}\Big[\frac{p^{2}c_0^{2}}{\epsilon}
|x|^{2p-2}+pc_0(N+p-2)|x|^{p-2}
\Big]\hat{\nu}_{\alpha}^2dx\Big\}\nonumber\\
\geq&\alpha^{\frac{2}{p}}\hat{\lambda}\ \ \text{as}\ \ \alpha\rightarrow\infty,\nonumber
\end{align}
where $\hat{\nu}_{\alpha}(x)=\alpha^{-\frac{N}{2p}}u_{\alpha}(\alpha^{-\frac{1}{p}}x)$, and $u_{\alpha}(x)\equiv 0$ in $\mathbb{R}^{N}\backslash\Omega$. This gives the lower  estimate of $\lambda(\alpha)$ as $\alpha\rightarrow\infty$.

We next discuss the upper estimate of $\lambda(\alpha)$ as $\alpha\rightarrow\infty$.
Choose a non-negative function $\varphi(x)\in C_{0}^{\infty}(\mathbb{R}^{N})$ such that $\varphi(x)=1$ for $|x|\leq\frac{r_{0}}{2}$ and $\varphi(x)=0$ for $|x|\geq r_{0}$, where $r_{0}>0$ is as in ($M$).
For any $\alpha>0$, set
\begin{equation}\label{2.13}
\phi_{\alpha}(x):=A_{\alpha}\alpha^{\frac{N}{2p}}\varphi(x)
\hat{u}(\alpha^{\frac{1}{p}}x),\ \ x\in\Omega,
\end{equation}
where $\hat{u}(x)>0$ denotes the unique positive minimizer of $\hat{\lambda}$ defined by (\ref{2.3}), and $A_{\alpha}>0$ is chosen such that $\int_{\Omega}|\phi_{\alpha}|^{2}dx=1$. By the exponential decay of Lemma \ref{lem:2.1}, we then have
\begin{equation}\label{2.14}
\frac{1}{A^{2}_{\alpha}}=1+O(\alpha^{-\infty})\ \ \text{as}\ \ \alpha\rightarrow\infty,
\end{equation}
where and below we use $f(t)=O(t^{-\infty})$ to denote a function $f$ satisfying $\lim_{t\rightarrow\infty}|f(t)|t^s=0$ for all $s>0$. Under the assumptions of Lemma \ref{lem:2.2},
direct calculations thus yield that
\begin{equation}\label{2.49}
\begin{split}
&\int_{\Omega}\epsilon|\nabla \phi_{\alpha}|^{2}dx+\int_{\Omega}\Big[\alpha p(N+p-2)g(x)|x|^{p-2}+\frac{\alpha^{2}p^{2}}{\epsilon}g^{2}(x)|x|^{2p-2}\Big]\phi^{2}_{\alpha}dx\\
=&\alpha^{\frac{2}{p}}\Big\{\int_{\mathbb{R}^{N}}\epsilon|\nabla \hat{u}|^{2}dx
+\int_{\mathbb{R}^{N}}\Big[p(N+p-2)g(\alpha^{-\frac{1}{p}}x)|x|^{p-2}\\
&+\frac{p^{2}}{\epsilon}g^{2}(\alpha^{-\frac{1}{p}}x)|x|^{2p-2}\Big]\hat{u}^{2}dx\Big\}
+O(\alpha^{-\infty})\\
=&\alpha^{\frac{2}{p}}\Big\{\int_{\mathbb{R}^{N}}\epsilon|\nabla \hat{u}|^{2}dx
+\int_{\mathbb{R}^{N}}\Big[pc_{0}(N+p-2)|x|^{p-2}
+\frac{p^{2}c_{0}^{2}}{\epsilon}|x|^{2p-2}\Big]\hat{u}^{2}dx+o(1)\Big\}\\
=&\alpha^{\frac{2}{p}}\big(\hat{\lambda}+o(1)\big)
\ \ \text{as}\ \ \alpha\rightarrow\infty,
\end{split}
\end{equation}
and
\begin{equation}\label{2.50}
\begin{split}
&\int_{\Omega} \Big\{\alpha\Big[|x|^{p}\Delta g+2p|x|^{p-2}(x\cdot\nabla g)\Big]+\frac{\alpha^{2}}{\epsilon}\Big[|x|^{2p}|\nabla g|^{2}+p|x|^{2p-2}\big(x\cdot\nabla g^{2}\big)\Big]\Big\}
\phi^{2}_{\alpha}dx\\
\leq&\alpha^{\frac{2}{p}}\int_{B_{r_{0}\alpha^{\frac{1}{2p}}}(0)} \Big\{|x|^{p}\big|\Delta g(\alpha^{-\frac{1}{p}}x)\big|+2p|x|^{p-1}\big|\nabla g(\alpha^{-\frac{1}{p}}x)\big|   \\
&\qquad +\frac{1}{\epsilon}\big[|x|^{2p}\big|\nabla g(\alpha^{-\frac{1}{p}}x)\big|^{2}+p|x|^{2p-1}\big|\nabla g^{2}(\alpha^{-\frac{1}{p}}x)\big|\big]\Big\}\hat{u}^{2}dx+O(\alpha^{-\infty})\\
\leq &C\alpha^{\frac{1-l}{p}}
\ \ \text{as}\ \ \alpha\rightarrow\infty,
\end{split}
\end{equation}
where $g(x)\equiv0$ in $\mathbb{R}^{N}\backslash\Omega$, $l>0$ and $r_{0}>0$ are as in ($M$), and $C>0$ is independent of $\alpha>0$.
Moreover, it follows from (\ref{2.9}), (\ref{2.13}) and (\ref{2.14}) that
\begin{equation}{\label{0.4}}
\int_{\Omega}V(x)\phi_{\alpha}^2dx=V(0)+o(1)\ \ \text{as}\ \ \alpha\rightarrow\infty.
\end{equation}
Applying (\ref{2.1}), we hence derive from (\ref{2.49})--(\ref{0.4}) that
\begin{equation}\label{2.15}
\begin{split}
\lambda(\alpha)\le &\int_{\Omega}\epsilon|\nabla \phi_{\alpha}|^2dx+\int_{\Omega}\Big[\frac{\alpha^{2}}{\epsilon}|\nabla m|^{2}
+\alpha\Delta m+V(x)\Big]\phi_{\alpha}^2dx\\
=&\int_{\Omega}\epsilon|\nabla \phi_{\alpha}|^2dx+\int_{\Omega}\Big\{\alpha\Big(2p(x\cdot\nabla g)+p(N+p-2)g(x)\Big)|x|^{p-2}\\
&+\frac{\alpha^{2}}{\epsilon}\Big[|x|^{2p}|\nabla g|^{2}+p|x|^{2p-2}(x\cdot\nabla g^{2})+p^{2}g^{2}(x)|x|^{2p-2}\Big]\\
&+\alpha|x|^p\Delta g+V(x)\Big\}\phi_{\alpha}^2dx\\
\leq&\alpha^{\frac{2}{p}}\big(\hat{\lambda}+o(1)\big)
\ \ \text{as}\ \ \alpha\rightarrow\infty.
\end{split}
\end{equation}
We thus conclude from (\ref{2.12}) and (\ref{2.15}) that (\ref{1.2}) holds true,
and the proof of Lemma \ref{lem:2.2} is therefore complete.\qed

Employing Lemmas \ref{lem:2.1} and \ref{lem:2.2}, we next derive the following limiting behavior of the unique principal eigenfunction $u_{\alpha}$ for $(\ref{1.1})$ as $\alpha\rightarrow\infty$, where we always define $u_\alpha(x)\equiv0$ for $x\in\mathbb{R}^{N}\backslash\Omega$.

\begin{prop}\label{prop2.3}
Suppose $0\leq V (x)\in C^{\gamma}(\bar{\Omega})~(0<\gamma<1)$, and assume $m(x)=g(x)|x|^p\geq0$, where $p\geq2$ and $g(x)$ satisfies ($M$) with $l>0$.
Then for any fixed $\epsilon>0$, the principal eigenfunction $u_{\alpha}$ of $(\ref{1.1})$  satisfies
\begin{equation}\label{1.3}
\begin{split}
\hat{w}_\alpha(x)&:=\alpha^{-\frac{N}{2p}}u_{\alpha}
(\alpha^{-\frac{1}{p}}x+x_\alpha)\\
&\rightarrow\hat{u}(x)>0
\ \ \text{strongly in}\ \ H^{1}(\mathbb{R}^{N})\cap L^{\infty}(\mathbb{R}^{N})
\ \ \text{as}\ \ \alpha\rightarrow\infty,
\end{split}
\end{equation}
where $\hat{u}$ is the unique positive minimizer of $\hat{\lambda}$ defined by (\ref{2.3}), $u_\alpha(x)\equiv0$ in $\mathbb{R}^{N}\backslash\Omega$,
and $x_\alpha\in\mathbb{R}^{N}$ is the unique global maximum point of $u_\alpha$ as $\alpha\rightarrow\infty$ and satisfies
\begin{equation}\label{1.8}
\lim_{\alpha\rightarrow\infty}\alpha^{\frac{1}{p}}|x_\alpha|=0.
\end{equation}
\end{prop}

\noindent\textbf{Proof.} For $\Omega_{\alpha}:=\big\{x\in\R^N:\, \alpha^{-\frac{1}{p}}x+x_\alpha\in\Omega\big\}$, set
\begin{equation}\label{2.17}
\hat{w}_\alpha(x):=\left\{
\begin{split}\alpha^{-\frac{N}{2p}}u_{\alpha}\big(\alpha^{-\frac{1}{p}}x+x_\alpha\big),\ \ &x\in\Omega_{\alpha};\\[2mm]
0, \qquad\qquad &x\in\mathbb{R}^{N}\backslash\Omega_{\alpha};
\end{split}
\right.
\end{equation}
where $x_\alpha$ is a global maximum point of $u_\alpha(x)$. It then follows from (\ref{1.6}) and (\ref{2.17}) that $\hat{w}_{\alpha}$ satisfies the following equation
\begin{align}
&-\epsilon\Delta \hat{w}_{\alpha}(x)+\Big[\big|x+\alpha^{\frac{1}{p}}x_\alpha\big|^{p} \Delta g(\alpha^{-\frac{1}{p}}x+x_\alpha) +
\frac{p^{2}}{\epsilon}g^{2}(\alpha^{-\frac{1}{p}}x+x_\alpha)\big|x+\alpha^{\frac{1}{p}}x_\alpha\big|^{2p-2}
\nonumber\\
&+p(N+p-2)g(\alpha^{-\frac{1}{p}}x+x_\alpha)\big|x+\alpha^{\frac{1}{p}}x_\alpha\big|^{p-2}
+\frac{1}{\epsilon}\big|x+\alpha^{\frac{1}{p}}x_\alpha\big|^{2p}\big|\nabla g(\alpha^{-\frac{1}{p}}x+x_\alpha)\big|^{2}\nonumber\\
&+2p\big|x+\alpha^{\frac{1}{p}}x_\alpha\big|^{p-2}
\big(x+\alpha^{\frac{1}{p}}x_\alpha\big)\cdot\nabla g(\alpha^{-\frac{1}{p}}x+x_\alpha)+\alpha^{-\frac{2}{p}}V\big(\alpha^{-\frac{1}{p}}x
+x_\alpha\big)
\label{2.26}\\
&+\frac{p}{\epsilon}\big|x+\alpha^{\frac{1}{p}}x_\alpha\big|^{2p-2}
\big(x+\alpha^{\frac{1}{p}}x_\alpha\big)\cdot\nabla g^{2}(\alpha^{-\frac{1}{p}}x+x_\alpha)
\Big]
\hat{w}_{\alpha}(x)\nonumber\\
&=\alpha^{-\frac{2}{p}}\lambda(\alpha)\hat{w}_{\alpha}(x)\ \ \,\text{in}\ \ \Omega_{\alpha}.\nonumber
\end{align}
Since the origin is a global maximum point of $\hat{w}_{\alpha}(x)$, we deduce from ($M$) and (\ref{2.26}) that
\begin{equation}\label{2.39}
0\leq pc_{0}(N+p-2)|\alpha^{\frac{1}{p}}x_\alpha|^{p-2}
+\frac{p^{2}c_{0}^{2}}{\epsilon}|\alpha^{\frac{1}{p}}x_\alpha|^{2p-2}
\leq\alpha^{-\frac{2}{p}}\lambda(\alpha).
\end{equation}
We then deduce from (\ref{1.2}) and (\ref{2.39}) that there exists a constant $\delta>0$ such that
\begin{equation}\label{2.40}
\limsup_{\alpha\rightarrow\infty}|\alpha^{\frac{1}{p}}x_\alpha|\leq\delta.
\end{equation}
Therefore, there exist a subsequence $\{\alpha_{n}\}$ of $\{\alpha\}$ and a point $y_0\in\mathbb{R}^{N}$ such that
\begin{equation}\label{2.41}
\lim_{n\rightarrow\infty}\alpha_{n}=\infty \ \ \text{and}\ \
\lim_{n\rightarrow\infty}\alpha_{n}^{\frac{1}{p}}x_{\alpha_{n}}=y_0.
\end{equation}
Since it yields from (\ref{2.41}) that $x_{\alpha_{n}}$ approaches to the origin as $n\to\infty$, which is an interior point of $\Omega$, we conclude that $\lim_{n\rightarrow\infty}\Omega_{\alpha_{n}}=\mathbb{R}^{N}$.

One can check that
\begin{equation*}
(1+t)^{p_{0}}\leq1+C(p_{0})(t+t^{p_{0}}) \ \ \text{for any}\ \ t\geq0\ \ \text{and} \ \ p_{0}\geq2,
\end{equation*}
which further implies that
\begin{equation}\label{2.45}
(a+b)^{p_{0}}\leq a^{p_{0}}+C(p_{0})\big(a^{p_{0}-1}b+b^{p_{0}}\big) \ \ \text{for any}\ \ a>0\ \ \text{and} \ \ b>0.
\end{equation}
It then follows from (\ref{2.40}) and (\ref{2.45}) that for sufficiently large $\alpha_{n}>0,$
\begin{equation}\label{2.43}
\begin{split}
\int_{B_{2\delta}^{c}(0)}|x|^{2p-2}\hat{w}_{\alpha_{n}}^{2}dx
=&\int_{B_{2\delta}^{c}(0)}\big(|x|-|\alpha_{n}^{\frac{1}{p}}
x_{\alpha_{n}}|+|\alpha_{n}^{\frac{1}{p}}
x_{\alpha_{n}}|\big)^{2p-2}\hat{w}_{\alpha_{n}}^{2}dx\\
\leq&\int_{B_{2\delta}^{c}(0)}\big(|x|-|\alpha_{n}^{\frac{1}{p}}
x_{\alpha_{n}}|\big)^{2p-2}\hat{w}_{\alpha_{n}}^{2}dx
+C(p,\delta)\int_{B_{2\delta}^{c}(0)}\hat{w}_{\alpha_{n}}^{2}dx\\
&+C(p,\delta)\int_{B_{2\delta}^{c}(0)}\big(|x|-|\alpha_{n}^{\frac{1}{p}}
x_{\alpha_{n}}|\big)^{2p-3}\hat{w}_{\alpha_{n}}^{2}dx\\
\leq&C(p,\delta)\Big(\int_{\mathbb{R}^{N}}|x+\alpha_{n}^{\frac{1}{p}}
x_{\alpha_{n}}|^{2p-2}\hat{w}_{\alpha_{n}}^{2}dx\Big)^{\frac{2p-3}{2p-2}}\\
&+\int_{\mathbb{R}^{N}}|x+\alpha_{n}^{\frac{1}{p}}
x_{\alpha_{n}}|^{2p-2}\hat{w}_{\alpha_{n}}^{2}dx
+C(p,\delta),
\end{split}
\end{equation}
where $\delta>0$ is as in (\ref{2.40}). Note from ($M$), (\ref{1.2}) and (\ref{2.26}) that
\begin{equation}\label{2.18}
\begin{split}
&\hat{\lambda}+1\geq\alpha_{n}^{-\frac{2}{p}}\lambda(\alpha_{n})\\
&\geq\epsilon\int_{\mathbb{R}^{N}}|\nabla \hat{w}_{\alpha_{n}}|^2dx+\frac{p^{2}c^{2}_{0}}{\epsilon}\int_{\mathbb{R}^{N}}
\big|x+\alpha_{n}^{\frac{1}{p}}x_{\alpha_{n}}\big|^{2p-2}
\hat{w}_{\alpha_{n}}^2dx\ \ \text{as}\ \ n\rightarrow\infty.
\end{split}
\end{equation}
We then deduce from (\ref{2.43}) and (\ref{2.18}) that the sequence $\{\hat{w}_{\alpha_{n}}\}$ is bounded uniformly in $\mathcal{H}$ defined by (\ref{2.5}). Since the embedding $\mathcal{H}\hookrightarrow L^{q}(\mathbb{R}^{N})$ is compact  (cf. \cite[Lemma 3.1]{Z}) for $2\leq q<2^{*}$, where $2^{*}=\frac{2N}{N-2}$ if $N\geq3$ and $2^{*}=\infty$ if $N=1,2$, there exist a subsequence, still denoted by $\{\hat{w}_{\alpha_{n}}\}$, of $\{\hat{w}_{\alpha_{n}}\}$ and a function $0\leq\hat{w}_{0}\in \mathcal{H}$ such that
\begin{equation}\label{0.7}
\begin{split}
\hat{w}_{\alpha_{n}}\rightharpoonup\hat{w}_{0}\ \ \text{weakly in}\ \ H^{1}(\mathbb{R}^{N})\ \ \text{as}\ \ n\rightarrow\infty,\qquad\qquad\\
\hat{w}_{\alpha_{n}}\rightarrow\hat{w}_{0} \ \ \text{strongly in}\ \ L^{q}(\mathbb{R}^{N})~(2\leq q<2^*)\ \ \text{as}\ \ n\rightarrow\infty,
\end{split}
\end{equation}
which also implies that $\int_{\mathbb{R}^{N}}|\hat{w}_{0}|^{2}dx=1$, due to the fact that $\int_{\mathbb{R}^{N}}|\hat{w}_{\alpha_{n}}|^{2}dx\equiv1$ holds for all $n>0$.

By Fatou's lemma, we now obtain from ($M$), (\ref{1.2}), (\ref{2.26}) and (\ref{2.41}) that
\begin{equation}\label{2.2}
\begin{split}
\hat{\lambda}=&\lim_{n\rightarrow\infty}
\alpha_{n}^{-\frac{2}{p}}\lambda(\alpha_{n})\\
=&\lim_{n\rightarrow\infty}\Big\{
\epsilon\int_{\mathbb{R}^{N}}\big|\nabla \hat{w}_{\alpha_{n}}\big|^2dx+\int_{\mathbb{R}^{N}}
\Big[\big|x+\alpha_{n}^{\frac{1}{p}}x_{\alpha_{n}}\big|^{p}\Delta g\big(\alpha^{-\frac{1}{p}}_{n}x+x_{\alpha_{n}}\big)\\
&+p(N+p-2)g\big(\alpha^{-\frac{1}{p}}_{n}x+x_{\alpha_{n}}\big)
\big|x+\alpha_{n}^{\frac{1}{p}}x_{\alpha_{n}}\big|^{p-2}+
\alpha_{n}^{-\frac{2}{p}}V\big(\alpha_{n}^{-\frac{1}{p}}x+x_{\alpha_{n}}\big)\\
&+\frac{1}{\epsilon}\big|x+\alpha_{n}^{\frac{1}{p}}x_{\alpha_{n}}\big|^{2p}\big|\nabla
g\big(\alpha^{-\frac{1}{p}}_{n}x+x_{\alpha_{n}}\big)\big|^{2}
+\frac{p^{2}}{\epsilon}g^{2}\big(\alpha^{-\frac{1}{p}}_{n}x+x_{\alpha_{n}}\big)
\big|x+\alpha_{n}^{\frac{1}{p}}x_{\alpha_{n}}\big|^{2p-2}\\
&+\frac{p}{\epsilon}\big|x+\alpha_{n}^{\frac{1}{p}}x_{\alpha_{n}}\big|^{2p-2}
\big(x+\alpha_{n}^{\frac{1}{p}}x_{\alpha_{n}}\big)\cdot\nabla
g^{2}\big(\alpha^{-\frac{1}{p}}_{n}x+x_{\alpha_{n}}\big)\\
&+2p\big|x+\alpha_{n}^{\frac{1}{p}}x_{\alpha_{n}}\big|^{p-2}
\big(x+\alpha_{n}^{\frac{1}{p}}x_{\alpha_{n}}\big)
\cdot\nabla
g\big(\alpha^{-\frac{1}{p}}_{n}x+x_{\alpha_{n}}\big)
\Big]
\hat{w}_{\alpha_{n}}^2(x)\chi_{\Omega_{\alpha_{n}}}(x)dx
\Big\}\\
\geq&\liminf_{n\rightarrow\infty}
\int_{\mathbb{R}^{N}}\epsilon\big|\nabla \hat{w}_{\alpha_{n}}\big|^2dx+\liminf_{n\rightarrow\infty}
\int_{\mathbb{R}^{N}}\Big[\frac{p^{2}}{\epsilon}g^{2}\big(\alpha^{-\frac{1}{p}}_{n}x+x_{\alpha_{n}}
\big)\big|x
+\alpha_{n}^{\frac{1}{p}}x_{\alpha_{n}}\big|^{2p-2}\\
&+p(N+p-2)g\big(\alpha^{-\frac{1}{p}}_{n}x+x_{\alpha_{n}}\big)\big|x+\alpha_{n}^{\frac{1}{p}}
x_{\alpha_{n}}\big|^{p-2}\Big]
\hat{w}_{\alpha_{n}}^2(x)\chi_{\Omega_{\alpha_{n}}}(x)dx,
\end{split}
\end{equation}
which then implies from (\ref{0.7}) that
\begin{equation}\label{2.21}
\begin{split}
\hat{\lambda}\geq& \epsilon\int_{\mathbb{R}^{N}}\big|\nabla \hat{w}_{0}\big|^2dx\\
&+\int_{\mathbb{R}^{N}}
\Big[\frac{p^{2}c_{0}^{2}}{\epsilon}\big|x
+y_0\big|^{2p-2}+pc_{0}(N+p-2)\big|x+y_0\big|^{p-2}\Big]\hat{w}_{0}^2(x)dx \\
=&\epsilon\int_{\mathbb{R}^{N}}|\nabla \hat{w}_{0}|^2dx+\int_{\mathbb{R}^{N}}\Big[\frac{p^{2}c_{0}^{2}}{\epsilon}|x|^{2p-2}
+pc_{0}(N+p-2)|x|^{p-2}\Big]\hat w_{0}^2(x-y_0)dx\\
\geq&\hat{\lambda}.
\end{split}\end{equation}
Here $\hat w_{0}\geq0$ is defined by (\ref{0.7}) and $y_0\in\mathbb{R}^{N}$ is as in (\ref{2.41}). This thus yields from Lemma \ref{lem:2.1} that $\hat w_0(x-y_0)\equiv \hat{u}(x)>0$ in $\mathbb{R}^{N}$, where $\hat{u}(x)$ is the unique positive minimizer of $\hat{\lambda}$ and $\hat{u}(0)=\max_{\mathbb{R}^{N}}\hat{u}(x)$. However, one can get from (\ref{2.17}) that the origin is always a maximum point of $\hat{w}_{\alpha_{n}}(x)$ for all $n>0$, which gives that the origin is a maximum point of $\hat{w}_{0}(x)$ in $\R^N$. We hence have $y_0=0$, $\hat w_0(x)\equiv \hat{u}(x)>0$ in $\mathbb{R}^{N}$, and it follows from (\ref{2.41}) that (\ref{1.8}) holds true.

We next prove that
\begin{equation}\label{2.22}
\hat{w}_{\alpha_{n}}(x)\rightarrow\hat{w}_{0}(x)=\hat{u}(x)\ \ \text{strongly in}\ \ H^{1}(\mathbb{R}^{N})\ \ \text{as}\ \ n\rightarrow\infty.
\end{equation}
Indeed, it follows from ($M$), (\ref{1.2}) and (\ref{2.26}) that
\begin{equation}\label{2.27}
-\Delta \hat{w}_{\alpha_{n}}(x)-\frac{2\hat{\lambda}}
{\epsilon}\hat{w}_{\alpha_{n}}(x)\leq0\ \ \text{in}\ \ \Omega_{\alpha_{n}}\ \ \text{as}\ \ n\rightarrow\infty.
\end{equation}
Since $\lim_{n\rightarrow\infty}\Omega_{\alpha_{n}}=\mathbb{R}^{N}$, we obtain that  for any $R>0$, we have $B_{R+2}(0)\subset\Omega_{\alpha_{n}}$ as $n\rightarrow\infty$.
Applying De Giorgi-Nash-Moser theory  (cf. \cite[Theorem 4.1]{Han}), we thus conclude from (\ref{2.27}) that there exists a constant $C>0$, independent of $n>0$, such that
\begin{equation}\label{2.28}
\max\limits_{B_{\frac{1}{2}}(\rho)}\hat{w}_{\alpha_{n}}(x)\leq C\Big(\int_{{B_{1}(\rho)}}|\hat{w}_{\alpha_{n}}(x)|^{2} dx\Big)^\frac{1}{{2}}\ \ \mbox{as}\ \ n\rightarrow\infty,
\end{equation}
where $\rho\in B_{R+1}(0)$.
Since $\int_{\mathbb{R}^{N}}|\hat{w}_{\alpha_{n}}|^{2}dx\equiv 1$ holds for all $n>0$, it follows from (\ref{2.28}) that
\begin{equation}\label{2.32}
\{\hat{w}_{\alpha_{n}}\}\ \ \text{is bounded uniformly in}\ \ H^{1}\big(B_{R+1}(0)\big)\cap L^{\infty}\big(B_{R+1}(0)\big)
\ \ \text{as}\ \ n\rightarrow\infty.
\end{equation}

On the other hand, we deduce from (\ref{1.2}), (\ref{2.26}) and (\ref{2.40}) that there exists a sufficiently large constant $R_{1}>0$, independent of $n>0$, such that
\begin{equation}\label{2.36}
-\Delta \hat{w}_{\alpha_{n}}+\hat{w}_{\alpha_{n}}\leq0\ \ \text{uniformly in}\ \ \Omega_{\alpha_{n}}\backslash B_{R_{1}}(0)\ \ \text{as}\ \ n\rightarrow\infty.
\end{equation}
Applying the comparison principle, one can deduce from (\ref{2.32}) and (\ref{2.36}) that there exists a constant $C>0$, independent of $n>0$, such that
\begin{equation}\label{2.37}
|\hat{w}_{\alpha_{n}}(x)|\leq Ce^{-|x|}\ \ \text{uniformly in}\ \ \Omega_{\alpha_{n}}\backslash B_{R_{1}}(0)\ \ \text{as}\ \ n\rightarrow\infty.
\end{equation}
Therefore, we derive from ($M$), (\ref{0.7}), (\ref{2.32}) and (\ref{2.37}) that
\begin{equation}\label{2.46}
\begin{split}
&\lim_{n\rightarrow\infty}\int_{\mathbb{R}^{N}}\Big[\frac{p^{2}}{\epsilon}g^{2}
(\alpha^{-\frac{1}{p}}_{n}x+x_{\alpha_{n}})\big|x
+\alpha_{n}^{\frac{1}{p}}x_{\alpha_{n}}\big|^{2p-2}\\
&+p(N+p-2)g(\alpha^{-\frac{1}{p}}_{n}x+x_{\alpha_{n}})
\big|x+\alpha_{n}^{\frac{1}{p}}x_{\alpha_{n}}\big|^{p-2}\Big]
\hat{w}_{\alpha_{n}}^2(x)\chi_{\Omega_{\alpha_{n}}}(x)dx\\
&=\int_{\mathbb{R}^{N}}\Big[\frac{p^{2}c_{0}^{2}}
{\epsilon}|x|^{2p-2}+pc_{0}(N+p-2)|x|^{p-2}\Big]\hat{w}_{0}^2(x)dx,
\end{split}
\end{equation}
and
\begin{equation}\label{2.51}
\begin{split}
&\lim_{n\rightarrow\infty}\int_{\mathbb{R}^{N}}
\Big[\big|x+\alpha_{n}^{\frac{1}{p}}x_{\alpha_{n}}\big|^{p}\Delta g(\alpha^{-\frac{1}{p}}_{n}x+x_{\alpha_{n}})
+\frac{1}{\epsilon}\big|x+\alpha_{n}^{\frac{1}{p}}x_{\alpha_{n}}\big|^{2p}\big|\nabla
g(\alpha^{-\frac{1}{p}}_{n}x+x_{\alpha_{n}})\big|^{2}\\
&+
\alpha_{n}^{-\frac{2}{p}}V\big(\alpha_{n}^{-\frac{1}{p}}x+x_{\alpha_{n}}\big)
+\frac{p}{\epsilon}\big|x+\alpha_{n}^{\frac{1}{p}}x_{\alpha_{n}}\big|^{2p-2}
\big(x+\alpha_{n}^{\frac{1}{p}}x_{\alpha_{n}}\big)\cdot\nabla
g^{2}(\alpha^{-\frac{1}{p}}_{n}x+x_{\alpha_{n}})\\
&+2p\big|x+\alpha_{n}^{\frac{1}{p}}x_{\alpha_{n}}\big|^{p-2}
\big(x+\alpha_{n}^{\frac{1}{p}}x_{\alpha_{n}}\big)\cdot\nabla
g(\alpha^{-\frac{1}{p}}_{n}x+x_{\alpha_{n}})
\Big]\hat{w}_{\alpha_{n}}^2(x)\chi_{\Omega_{\alpha_{n}}}(x)dx=0.
\end{split}
\end{equation}
We then deduce from (\ref{2.2}), (\ref{2.21}), (\ref{2.46}) and (\ref{2.51}) that
$$\lim_{n\rightarrow\infty}\int_{\mathbb{R}^{N}}|\nabla \hat{w}_{\alpha_{n}}|^{2}dx
=\int_{\mathbb{R}^{N}}|\nabla \hat{w}_{0}|^{2}dx=\int_{\mathbb{R}^{N}}|\nabla \hat{u}|^{2}dx,$$
which further implies that (\ref{2.22}) holds true.

We now claim that
\begin{equation}\label{2.23}
\hat{w}_{\alpha_{n}}(x)\rightarrow\hat{u}(x)\ \ \text{strongly in}\ \ L^{\infty}(\mathbb{R}^{N})\ \ \text{as}\ \ n\rightarrow\infty.
\end{equation}
In fact, for the case $N=1$, the above claim is trivial in view of (\ref{2.22}), due to the fact that the embedding $H^1(\mathbb{R})\hookrightarrow L^\infty(\mathbb{R})$ is continuous (cf. \cite[Theorem 8.5]{Lieb1}).
We next consider the case $N\geq2$.
By the exponential decay (\ref{2.9}) and (\ref{2.37}), to prove the claim (\ref{2.23}), the rest is to show that
\begin{equation}\label{2.38}
\hat{w}_{\alpha_{n}}(x)\rightarrow\hat{u}(x)\ \ \text{strongly in}\ \ L_{loc}^{\infty}(\mathbb{R}^{N})\ \ \text{as}\ \ n\rightarrow\infty.
\end{equation}
Indeed, setting
\begin{equation*}
\begin{split}
\hat{G}_{\alpha_{n}}(x):=&\frac{\alpha_{n}^{-\frac{2}{p}}
\lambda(\alpha_{n})}{\epsilon}\hat{w}_{\alpha_{n}}(x)-
\frac{1}{\epsilon}\Big\{2p\big|x+\alpha_{n}^{\frac{1}{p}}x_{\alpha_{n}}\big|^{p-2}
\big(x+\alpha_{n}^{\frac{1}{p}}x_{\alpha_{n}}\big)\cdot\nabla g(\alpha_{n}^{-\frac{1}{p}}x+x_{\alpha_{n}})\\
&+\frac{1}{\epsilon}\big|x+\alpha_{n}^{\frac{1}{p}}x_{\alpha_{n}}\big|^{2p}\big|\nabla g(\alpha_{n}^{-\frac{1}{p}}x+x_{\alpha_{n}})\big|^{2}+
\frac{p^{2}}{\epsilon}g^{2}(\alpha_{n}^{-\frac{1}{p}}x+x_{\alpha_{n}})
\big|x+\alpha_{n}^{\frac{1}{p}}x_{\alpha_{n}}\big|^{2p-2}\\
&+p(N+p-2)g(\alpha_{n}^{-\frac{1}{p}}x+x_{\alpha_{n}})\big|x+\alpha_{n}^{\frac{1}{p}}
x_{\alpha_{n}}\big|^{p-2}
+\alpha_{n}^{-\frac{2}{p}}V\big(\alpha_{n}^{-\frac{1}{p}}x
+x_{\alpha_{n}}\big)\\
&+\frac{p}{\epsilon}\big|x+\alpha_{n}^{\frac{1}{p}}x_{\alpha_{n}}\big|^{2p-2}
\big(x+\alpha_{n}^{\frac{1}{p}}x_{\alpha_{n}}\big)\cdot\nabla g^{2}(\alpha_{n}^{-\frac{1}{p}}x+x_{\alpha_{n}})\\
&+\big|x+\alpha_{n}^{\frac{1}{p}}x_{\alpha_{n}}\big|^{p}\Delta g(\alpha_{n}^{-\frac{1}{p}}x+x_{\alpha_{n}})
\Big\}\hat{w}_{\alpha_{n}}(x),
\end{split}
\end{equation*}
we get from (\ref{2.26}) that
\begin{equation}\label{2.31}
-\Delta\hat{w}_{\alpha_{n}}(x)=\hat{G}_{\alpha_{n}}(x) \ \ \mbox{in}\ \ \Omega_{\alpha_{n}}.
\end{equation}
For any $R>0$ and $r>N$, we derive from (\ref{1.2}), (\ref{2.40}) and (\ref{2.32}) that the sequence $\{\hat{G}_{\alpha_{n}}\}$ is bounded uniformly in $L^{r}(B_{R+1}(0))$ as $n\rightarrow\infty$. It then follows from (\ref{2.31}) and
\cite[Theorem 9.11]{GT} that
the sequence
$\{\hat{w}_{\alpha_{n}}\}$ is also bounded uniformly in $W^{2,r}(B_{R}(0))$ as $n\rightarrow\infty$.
Since $r>N$, we deduce from \cite[Theorem 7.26]{GT} that the sequence
$\{\hat{w}_{\alpha_{n}}\}$ is bounded uniformly in $C^{1}(B_{R}(0))$ as $n\rightarrow\infty$, and
there exists a subsequence, still denoted by $\{\hat{w}_{\alpha_{n}}\}$, of $\{\hat{w}_{\alpha_{n}}\}$ such that for some $w(x)\in C^{1}(B_{R}(0))$,
\begin{equation}\label{2.34}
\hat{w}_{\alpha_{n}}(x)\rightarrow w(x) \ \ \mbox{strongly in}\ \ C^{1}(B_{R}(0))\ \ \mbox{as}\ \  n\rightarrow\infty.
\end{equation}
Since $R>0$ is arbitrary, we now conclude from (\ref{2.22}) and (\ref{2.34})  that (\ref{2.38}) holds true. This further implies that the claim (\ref{2.23}) is true. Together with (\ref{2.22}), this hence proves (\ref{1.3}).

We finally prove the uniqueness of the global maximum point $x_{\alpha_{n}}$ for $u_{\alpha_{n}}$ as $n\rightarrow\infty$. Because it yields from (\ref{2.34}) that the sequence $\{\hat{w}_{\alpha_{n}}\}$ is bounded uniformly in $C_{loc}^{1}(\mathbb{R}^{N})$ as $n\rightarrow\infty$, it follows from the Schauder estimate \cite[Theorem 6.2]{GT} that the sequence $\{\hat{w}_{\alpha_{n}}\}$ is also bounded uniformly in $C_{loc}^{2,\varsigma}(\mathbb{R}^{N})$ as $n\rightarrow\infty$ for some $0<\varsigma<1$. Therefore, up to a subsequence if necessary, there exists a function $\bar{w}(x)\in C_{loc}^{2}(\mathbb{R}^{N})$ such that
$$\hat{w}_{\alpha_{n}}(x)\rightarrow\bar{w}(x)\ \ \text{in}\ \ C_{loc}^{2}(\mathbb{R}^{N})\ \ \text{as}\ \ n\rightarrow\infty,$$
and we further deduce from (\ref{2.23}) that
\begin{equation}\label{2.7}
\hat{w}_{\alpha_{n}}(x)\rightarrow\bar{w}(x)\equiv\hat{u}(x)\ \ \text{in}\ \ C_{loc}^{2}(\mathbb{R}^{N})\ \ \text{as}\ \ n\rightarrow\infty.
\end{equation}
Since the origin is the unique global maximum point of $\hat{u}(x)$, we derive from (\ref{2.9}), (\ref{2.37}) and (\ref{2.7}) that all global maximum points of $\hat{w}_{\alpha_{n}}$ must stay in a small ball $B_{\delta_{0}}(0)$ as $n\rightarrow\infty$ for some small $\delta_{0}>0$. It also yields from $\hat{u}''(0)<0$  that $\hat{u}''(r)<0$ holds for all $0<r<\delta_{0}$. We then reduce from \cite[Lemma 4.2]{Ni} that $\hat{w}_{\alpha_{n}}$ has a unique maximum point for sufficiently large $n>0$, which is just the origin. This thus proves the uniqueness of $x_{\alpha_{n}}$ as $n\rightarrow\infty$.

Because the above results are independent of the subsequence $\{\hat{w}_{\alpha_{n}}\}$ that we choose, we conclude  that both (\ref{1.3}) and (\ref{1.8}) hold essentially true for the whole sequence $\{\hat{w}_{\alpha}\}$. The proof of Proposition \ref{prop2.3} is therefore complete.\qed

\section{Refined Limiting Profiles of  $(\lambda(\alpha), u_{\alpha})$ as $\alpha\to\infty$}

The purpose of this section is to prove Theorem \ref{cor 1.2} on the refined limiting profiles of the principal eigenpair $(\lambda(\alpha),u_{\alpha})$ for (\ref{1.1}) as $\alpha\rightarrow\infty$,  no matter whether the minimum points of $V(x)$ and $m(x)$  coincide with each other or not.
We begin with the following second-order asymptotic expansion of $\lambda(\alpha)$ as $\alpha\rightarrow\infty$.

\begin{lem}\label{lem3.1}
Suppose $0\leq V (x)\in C^{\gamma}(\bar{\Omega})~(0<\gamma<1)$, and assume  $m(x)=g(x)|x|^p\geq0$, where $p\in\{2\}\cup(3,\infty)$ and $g(x)$ satisfies ($M$) for some  $l>1$.
Then for any fixed $\epsilon>0$, the principal eigenvalue $\lambda(\alpha)$ of (\ref{1.1}) satisfies
\begin{equation}\label{3.2}
\lambda(\alpha)=\alpha^{\frac{2}{p}}\hat{\lambda}+V(0)+o(1)\ \ \text{as}\ \ \alpha\rightarrow\infty,
\end{equation}
where $\hat{\lambda}>0$ is given by (\ref{2.3}).
\end{lem}

\noindent\textbf{Proof.}
Define
\begin{equation}\label{3.42}
\tilde{w}_\alpha(x):=\left\{
\begin{split}\alpha^{-\frac{N}{2p}}u_{\alpha}
\big(\alpha^{-\frac{1}{p}}x\big)>0,\ \ &x\in\Omega'_{\alpha}:=\big\{x\in\R^N:\,\alpha^{-\frac{1}{p}}x\in\Omega\big\};\\[2mm]
0, \qquad\qquad &x\in\mathbb{R}^{N}\backslash\Omega'_{\alpha};
\end{split}\right.
\end{equation}
where $u_{\alpha}>0$ is the unique positive principal eigenfunction of (\ref{1.1}).
It then follows from (\ref{1.6}) and (\ref{3.42}) that $\tilde{w}_{\alpha}$ satisfies
\begin{equation}\label{3.1}
\begin{split}
&-\epsilon\Delta \tilde{w}_{\alpha}+\Big\{|x|^{p}\Delta g(\alpha^{-\frac{1}{p}}x)+\big[2p\big(x\cdot\nabla g(\alpha^{-\frac{1}{p}}x)\big)+p(N+p-2)g(\alpha^{-\frac{1}{p}}x)\big]|x|^{p-2}
\\
&+\frac{1}{\epsilon}|x|^{2p}|\nabla g(\alpha^{-\frac{1}{p}}x)|^{2}+\frac{p}{\epsilon}|x|^{2p-2}
\Big(x\cdot\nabla g^{2}(\alpha^{-\frac{1}{p}}x)\Big)+
\frac{p^{2}}{\epsilon}g^{2}(\alpha^{-\frac{1}{p}}x)|x|^{2p-2}
\\
&+\alpha^{-\frac{2}{p}}V\big(\alpha^{-\frac{1}{p}}x\big)\Big\}
\tilde{w}_{\alpha}=\alpha^{-\frac{2}{p}}\lambda(\alpha)\tilde{w}_{\alpha}
\ \ \text{in}\ \ \Omega'_{\alpha}.
\end{split}
\end{equation}
Because the origin is an interior point of $\Omega$, we get that $\lim_{\alpha\rightarrow\infty}\Omega'_{\alpha}=\mathbb{R}^{N}$.
Since it yields from Proposition \ref{prop2.3} that $\lim_{\alpha\rightarrow\infty}\alpha^{\frac{1}{p}}|x_\alpha|=0$, one can deduce from (\ref{1.3}) that
\begin{equation}\label{3.47}
\tilde{w}_\alpha(x)\rightarrow\hat{u}(x)
\ \ \text{strongly in}\   L^{\infty}(\mathbb{R}^{N})\ \ \text{as}\ \ \alpha\rightarrow\infty.
\end{equation}
Similar to (\ref{2.37}), we also derive from (\ref{3.1}) that
there exists a constant $C>0$, independent of $\alpha>0$, such that
\begin{equation}\label{3.48}
|\tilde{w}_{\alpha}(x)|\leq Ce^{-|x|}\ \ \text{uniformly in}\ \ \Omega'_{\alpha}\backslash B_{R}(0)\ \ \text{as}\ \ \alpha\rightarrow\infty,
\end{equation}
where $R>0$ is large enough, and $\Omega'_{\alpha}$ is as in (\ref{3.42}). Furthermore,
applying the local elliptic estimate (cf. \cite[(3.15)]{GT}) to (\ref{3.1}), we deduce from (\ref{1.2}) and (\ref{3.48}) that there exists a constant $C>0$, independent of $\alpha>0$, such that
\begin{equation}\label{3.13M}
|\nabla\tilde{w}_{\alpha}(x)|\leq Ce^{-\frac{|x|}{2}}\ \ \text{uniformly in}\ \ \Omega'_{\alpha}\backslash B_{R}(0)\ \ \text{as}\ \ \alpha\rightarrow\infty,
\end{equation}
where $R>0$ is large enough.

Let $\hat{u}(x)>0$ be the unique positive minimizer of (\ref{2.3}).
Define
\begin{equation}\label{4}
\mathcal{L}:=-\epsilon\Delta+\Big[\frac{p^{2}c_{0}^{2}}{\epsilon}
|x|^{2p-2}+pc_{0}(N+p-2)|x|^{p-2}-\hat{\lambda}\Big]
\ \ \text{in}\  \ \mathbb{R}^{N},
\end{equation}
where $c_{0}>0$ is as in ($M$).
It then follows from (\ref{2.10}) that
\begin{equation}\label{0.2}
\mathcal{L}\hat{u}=0\ \ \text{in}\ \ \mathbb{R}^{N},\ \ \nabla\hat{u}(0)=0.
\end{equation}
We now claim that
\begin{equation}\label{A2}
ker\mathcal{L}=span\{\hat{u}\}\ \ \ \text{and}\ \ \ \langle\mathcal{L}u, u\rangle\geq0\ \ \ \text{for all}\ \  u\in\mathcal{H},
\end{equation}
where the space $\mathcal{H}$ is as in (\ref{2.5}).
In fact, suppose there exists a nonzero function $u_0\in\mathcal{H}$ such that
$\mathcal{L}u_0=0$ in $\mathbb{R}^{N}$. We then have
$$\hat{\lambda}=\frac{1}{\|u_0\|_{L^{2}(\mathbb{R}^{N})}^{2}}
\Big\{\int_{\mathbb{R}^{N}}\epsilon|\nabla u_0|^2dx
+\int_{\mathbb{R}^{N}}\Big[\frac{p^{2}c_{0}^{2}}
{\epsilon}|x|^{2p-2}+pc_{0}(N+p-2)|x|^{p-2}\Big]u_0^2dx\Big\}. $$
It thus follows from Lemma
\ref{lem:2.1} $(1)$ that there exists a constant $c_{1}\in\mathbb{R}$ such that
 $u_0=c_{1}\hat{u}$. Applying (\ref{0.2}),
we thus conclude that $ker\mathcal{L}=span\{\hat{u}\}$. On the other hand, if there exists $u_1\in\mathcal{H}$ such that
$\langle\mathcal{L}u_1, u_1\rangle<0$, then we have
$$\hat{\lambda}>\frac{1}{\|u_1\|_{L^{2}(\mathbb{R}^{N})}^{2}}
\Big\{\int_{\mathbb{R}^{N}}\epsilon|\nabla u_1|^2dx
+\int_{\mathbb{R}^{N}}\Big[\frac{p^{2}c_{0}^{2}}
{\epsilon}|x|^{2p-2}+pc_{0}(N+p-2)|x|^{p-2}\Big]u_1^2dx\Big\}. $$
This is a contradiction with the definition of $\hat{\lambda}$ in (\ref{2.3}). Therefore, the claim (\ref{A2}) holds true.

Inspired by (\ref{3.47}), we set
\begin{equation}\label{3.7M}
\tilde{v}_\alpha(x):=\tilde{w}_\alpha(x)-\hat{u}(x),\ \ x\in\mathbb{R}^{N}.
\end{equation}
It follows from (\ref{3.1}), (\ref{4}), (\ref{0.2}) and (\ref{3.7M}) that $\tilde{v}_\alpha(x)$ satisfies
\begin{equation}\label{3.12}
\begin{split}
\mathcal{L}\tilde{v}_{\alpha}
=&\Big\{-|x|^{p}\Delta g(\alpha^{-\frac{1}{p}}x)-2p|x|^{p-2}\big[x\cdot\nabla g(\alpha^{-\frac{1}{p}}x)\big]-\frac{1}{\epsilon}|x|^{2p}|\nabla g(\alpha^{-\frac{1}{p}}x)|^{2}\\
&+p(N+p-2)|x|^{p-2}\big[c_0-g(\alpha^{-\frac{1}{p}}x)\big]
-\frac{p}{\epsilon}|x|^{2p-2}
\big[x\cdot\nabla g^{2}(\alpha^{-\frac{1}{p}}x)\big]\\
&+\frac{p^{2}}{\epsilon}|x|^{2p-2}
\big[c_{0}^{2}-g^{2}(\alpha^{-\frac{1}{p}}x)\big]-\alpha^{-\frac{2}{p}}V(\alpha^{-\frac{1}{p}}x)
+\big[\alpha^{-\frac{2}{p}}\lambda(\alpha)-\hat{\lambda}\big]\Big\}
\tilde{w}_{\alpha}\\
:=&B_{\alpha}(x)\tilde{w}_{\alpha}(x)
\ \ \text{in}\ \ \Omega'_{\alpha},
\end{split}
\end{equation}
where $\Omega'_{\alpha}$ is as in (\ref{3.42}).
Recall from ($M$) that there exists a constant $C>0$, independent of $\alpha>0$, such that for $r_{0}>0$ and $l>1$ given by ($M$),
\begin{equation}\label{0.8}
|\Delta g(\alpha^{-\frac{1}{p}}x)|\leq C\alpha^{-\frac{l+2}{p}}|x|^{l}
\qquad \text{in}\ \ B_{r_{0}\alpha^{\frac{1}{2p}}}(0)\  \ \text{as}\ \ \alpha\rightarrow\infty,
\end{equation}
and
\begin{equation}\label{0.9}
|\nabla g(\alpha^{-\frac{1}{p}}x)|\leq C\alpha^{-\frac{l+1}{p}}|x|^{l}
\qquad  \text{in}\ \ B_{r_{0}\alpha^{\frac{1}{2p}}}(0)\  \ \text{as}\ \ \alpha\rightarrow\infty,
\end{equation}
together with
\begin{equation}\label{0.10}
\big|g(\alpha^{-\frac{1}{p}}x)-c_{0}\big|\leq C\alpha^{-\frac{l+1}{p}}|x|^{l+1}
\qquad  \text{in}\ \ B_{r_{0}\alpha^{\frac{1}{2p}}}(0)\  \ \text{as}\ \ \alpha\rightarrow\infty.
\end{equation}

We then derive from (\ref{2.9}), (\ref{3.48}) and (\ref{0.8})--(\ref{0.10}) that there exists a constant $C>0$, independent of $\alpha>0$, such that
\begin{equation}\label{3.13}
\begin{split}
0\leq&\Big|\int_{B_{r_{0}\alpha^{\frac{1}{2p}}}(0)}\Big\{-|x|^{p}\Delta g(\alpha^{-\frac{1}{p}}x)-2p|x|^{p-2}\big[x\cdot\nabla g(\alpha^{-\frac{1}{p}}x)\big]\\
&+p(N+p-2)|x|^{p-2}\big[c_0-g(\alpha^{-\frac{1}{p}}x)\big]
-\frac{p}{\epsilon}|x|^{2p-2}
\big[x\cdot\nabla g^{2}(\alpha^{-\frac{1}{p}}x)\big]\\
&-\frac{1}{\epsilon}|x|^{2p}|\nabla g(\alpha^{-\frac{1}{p}}x)|^{2}+\frac{p^{2}}{\epsilon}|x|^{2p-2}
\big[c_{0}^{2}-g^{2}(\alpha^{-\frac{1}{p}}x)\big]\Big\}
\tilde{w}_{\alpha}\hat{u}dx\Big|\\
\leq&C\alpha^{-\frac{l+1}{p}}\int_{B_{r_{0}\alpha^{\frac{1}{2p}}}(0)}
\Big[\alpha^{-\frac{1}{p}}|x|^{p+l}+2p|x|^{p+l-1}+p(N+p-2)|x|^{p+l-1}\\
&\qquad\qquad+\frac{p}{\epsilon}|x|^{2p+l-1}+\frac{1}{\epsilon}\alpha^{-\frac{l+1}{p}}|x|^{2p+2l}
+\frac{p^{2}}{\epsilon}|x|^{2p+l-1}\Big]
\tilde{w}_{\alpha}\hat{u}dx\\
\leq&C\alpha^{-\frac{l+1}{p}}\ \ \text{as}\ \ \alpha\rightarrow\infty,
\end{split}
\end{equation}
where $r_{0}>0$ is as in ($M$).
Multiplying (\ref{3.12}) by $\hat{u}(x)$ and integrating over $B_{r_{0}\alpha^{\frac{1}{2p}}}(0)\subseteq\Omega'_{\alpha}$, where $r_{0}>0$ is as in ($M$) and $\Omega'_{\alpha}$ is as in (\ref{3.42}), we further calculate from (\ref{2.9}), (\ref{3.48}), (\ref{3.13M}), (\ref{0.2}), (\ref{3.12}) and (\ref{3.13}) that for $l>1$ given by ($M$),
\begin{equation}\label{3.15}
\begin{split}
&\int_{B_{r_{0}\alpha^{\frac{1}{2p}}}(0)}B_{\alpha}(x)\tilde{w}_{\alpha}\hat{u}dx\\
=&\int_{B_{r_{0}\alpha^{\frac{1}{2p}}}(0)}\Big[-\alpha^{-\frac{2}{p}}
V(\alpha^{-\frac{1}{p}}x)
+\big(\alpha^{-\frac{2}{p}}\lambda(\alpha)-\hat{\lambda}\big)
\Big]\tilde{w}_{\alpha}\hat{u}dx+O(\alpha^{-\frac{l+1}{p}})\\
=&-\alpha^{-\frac{2}{p}}V(0)+o(\alpha^{-\frac{2}{p}})+
\big(\alpha^{-\frac{2}{p}}\lambda(\alpha)
-\hat{\lambda}\big)\big[1+o(1)\big]\ \ \text{as}\ \ \alpha\rightarrow\infty,
\end{split}
\end{equation}
and
\begin{equation}\label{3.14}
\begin{split}
&\int_{B_{r_{0}\alpha^{\frac{1}{2p}}}(0)}(\mathcal{L}\tilde{v}_{\alpha})\hat{u}dx\\
=&\int_{B_{r_{0}\alpha^{\frac{1}{2p}}}(0)}\tilde{v}_{\alpha}(\mathcal{L}\hat{u})dx
+\epsilon\int_{\partial B_{r_{0}\alpha^{\frac{1}{2p}}}(0)}
\Big(\tilde{v}_{\alpha}\frac{\partial\hat{u}}{\partial\nu}-\hat{u}
\frac{\partial\tilde{v}_{\alpha}}{\partial\nu}\Big)dS\\
=&o(e^{-r_{0}\alpha^{\frac{1}{2p}}})\ \ \text{as}\ \ \alpha\to\infty,
\end{split}
\end{equation}
where $\nu=(\nu^{1},\cdots,\nu^{N})$ denotes the outward unit normal vector of $\partial B_{r_{0}\alpha^{\frac{1}{2p}}}(0)$.

We hence obtain from (\ref{3.12}), (\ref{3.15}) and (\ref{3.14}) that
\begin{equation}\label{3.3}
-\alpha^{-\frac{2}{p}}\big[V(0)+o(1)\big]
+\big(\alpha^{-\frac{2}{p}}\lambda(\alpha)-\hat{\lambda}\big)
\big[1+o(1)\big]
=o(e^{-r_{0}\alpha^{\frac{1}{2p}}})
\ \ \text{as}\ \ \alpha\rightarrow\infty.
\end{equation}
This further implies that
\begin{equation}\label{3.17}
\alpha^{-\frac{2}{p}}\lambda(\alpha)-\hat{\lambda}
=\alpha^{-\frac{2}{p}}\big[V(0)+o(1)\big]
\ \ \text{as}\ \ \alpha\rightarrow\infty,
\end{equation}
and thus (\ref{3.2}) holds true. The proof of Lemma \ref{lem3.1} is therefore complete.
\qed

By applying Lemma \ref{lem3.1}, we next establish the following rough asymptotic expansion of the principal eigenfunction $u_{\alpha}$ as $\alpha\rightarrow\infty$.

\begin{lem}\label{lem3.2}
Under the assumptions of Lemma \ref{lem3.1}, the principal eigenfunction $u_{\alpha}$ of (\ref{1.1}) satisfies
\begin{equation}\label{3.201}
\begin{split}
\tilde{w}_{\alpha}(x):&=\alpha^{-\frac{N}{2p}}u_{\alpha}
(\alpha^{-\frac{1}{p}}x)\\
&=\Big(\int_{\mathbb{R}^{N}}\tilde{w}_{\alpha}\hat{u}dx\Big)
\hat{u}(x)+o(\alpha^{-\frac{2}{p}})\ \ \text{in}\ \ \mathbb{R}^{N}\ \ \text{as}\ \ \alpha\rightarrow\infty,
\end{split}
\end{equation}
where $u_\alpha(x)\equiv0$ in $\mathbb{R}^{N}\backslash\Omega$, and $\hat{u}(x)>0$ is the unique principal eigenfunction of  $\hat{\lambda}$ defined by (\ref{2.3}).
\end{lem}

\noindent\textbf{Proof.}
Define the operator $\mathcal{L}_{\alpha}$ by
\begin{equation}\label{3.58M}
\begin{split}
\mathcal{L}_{\alpha}:=&-\epsilon\Delta+\Big[\frac{p^{2}}{\epsilon}g^{2}(\alpha^{-\frac{1}{p}}x)
|x|^{2p-2}+p(N+p-2)g(\alpha^{-\frac{1}{p}}x)|x|^{p-2}\\
&+
\alpha^{-\frac{2}{p}}V(\alpha^{-\frac{1}{p}}x)-\alpha^{-\frac{2}{p}}\lambda(\alpha)\Big]
\ \ \text{in}\ \ \Omega'_{\alpha},
\end{split}
\end{equation}
where $\Omega'_{\alpha}$ is as in (\ref{3.42}).
We also define
\begin{equation}\label{3.58}
\eta_{\alpha}(x):=\tilde{w}_{\alpha}(x)
-\Big(\int_{\mathbb{R}^{N}}\tilde{w}_{\alpha}\hat{u}dx\Big)\hat{u}(x),\ \   x\in\mathbb{R}^{N},
\end{equation}
so that
\begin{equation}\label{3.58M1}
\int_{\mathbb{R}^{N}}\eta_{\alpha}\hat{u}dx\equiv0\ \ \text{for any}\ \, \alpha>0,
\end{equation}
where $\tilde{w}_{\alpha}(x)\geq0$ is as in (\ref{3.201}).
Recall from (\ref{3.47}) that
$$\int_{\mathbb{R}^{N}}\tilde{w}_{\alpha}\hat{u}dx\rightarrow1\ \ \text{as}\ \ \alpha\rightarrow\infty.$$

It follows from (\ref{2.10}), (\ref{3.1}) and (\ref{3.58M}) that
\begin{equation}\label{0.5}
\begin{split}
\mathcal{L}_{\alpha}\tilde{w}_{\alpha}=&-\Big[|x|^{p}\Delta g(\alpha^{-\frac{1}{p}}x)+\frac{p}{\epsilon}|x|^{2p-2}
\Big(x\cdot\nabla g^{2}(\alpha^{-\frac{1}{p}}x)\Big)
\\
&+\frac{1}{\epsilon}|x|^{2p}|\nabla g(\alpha^{-\frac{1}{p}}x)|^{2}+2p|x|^{p-2}\big(x\cdot\nabla g(\alpha^{-\frac{1}{p}}x)\big)\Big]\tilde{w}_{\alpha}\ \ \text{in}\ \ \Omega'_{\alpha},
\end{split}
\end{equation}
and
\begin{equation}\label{0.6}
\begin{split}
\mathcal{L}_{\alpha}\hat{u}=&\Big[
\alpha^{-\frac{2}{p}}V(\alpha^{-\frac{1}{p}}x)
+\frac{p^{2}}{\epsilon}|x|^{2p-2}\big(g^{2}(\alpha^{-\frac{1}{p}}x)-c^{2}_0\big)
-\big(\alpha^{-\frac{2}{p}}
\lambda(\alpha)-\hat{\lambda}\big)\\
&+p(N+p-2)|x|^{p-2}\big(g(\alpha^{-\frac{1}{p}}x)-c_0\big)
\Big]\hat{u}\ \ \text{in}\ \ \Omega'_{\alpha}.
\end{split}
\end{equation}
Additionally, it yields from Lemma \ref{lem3.1} that
\begin{equation}\label{0.11}
\begin{split}
&\alpha^{-\frac{2}{p}}V(\alpha^{-\frac{1}{p}}x)-\big(\alpha^{-\frac{2}{p}}
\lambda(\alpha)-\hat{\lambda}\big)\\
=&\alpha^{-\frac{2}{p}}\big[V(0)+o(1)\big]-\big(\alpha^{-\frac{2}{p}}
\lambda(\alpha)-\hat{\lambda}\big)\\
=&o(\alpha^{-\frac{2}{p}})
\ \ \text{in}\ \ B_{r_{0}\alpha^{\frac{1}{2p}}}(0)\  \ \text{as}\ \ \alpha\rightarrow\infty,
\end{split}
\end{equation}
where $r_{0}>0$ is as in ($M$).
We then calculate from (\ref{0.8})--(\ref{0.10}) and  (\ref{0.5})--(\ref{0.11}) that there exists a constant $C>0$, independent of $\alpha>0$, such that for $l>1$ given by ($M$),
\begin{equation}\label{3.59}
\begin{split}
|\mathcal{L}_{\alpha}\eta_{\alpha}|
=&\Big|\mathcal{L}_{\alpha}\tilde{w}_{\alpha}
-\Big(\int_{\mathbb{R}^{N}}\tilde{w}_{\alpha}\hat{u}dx\Big)
\mathcal{L}_{\alpha}\hat{u}\Big|\\
=&\Big|-\Big\{|x|^{p}\Delta g(\alpha^{-\frac{1}{p}}x)+2p|x|^{p-2}\big(x\cdot\nabla g(\alpha^{-\frac{1}{p}}x)\big)\\
&+\frac{1}{\epsilon}|x|^{2p}|\nabla g(\alpha^{-\frac{1}{p}}x)|^{2}+\frac{p}{\epsilon}|x|^{2p-2}
\Big(x\cdot\nabla g^{2}(\alpha^{-\frac{1}{p}}x)\Big)\Big\}\tilde{w}_{\alpha}
\\
&-\Big(\int_{\mathbb{R}^{N}}\tilde{w}_{\alpha}\hat{u}dx\Big)
\Big\{
\alpha^{-\frac{2}{p}}V(\alpha^{-\frac{1}{p}}x)
+\frac{p^{2}}{\epsilon}|x|^{2p-2}\big(g^{2}(\alpha^{-\frac{1}{p}}x)-c^{2}_0\big)\\
&+p(N+p-2)|x|^{p-2}\big(g(\alpha^{-\frac{1}{p}}x)-c_0\big)
-\big(\alpha^{-\frac{2}{p}}
\lambda(\alpha)-\hat{\lambda}\big)
\Big\}\hat{u}\Big|\\
\leq& \Big(\int_{\mathbb{R}^{N}}\tilde{w}_{\alpha}\hat{u}dx\Big)
\Big|\alpha^{-\frac{2}{p}}V(\alpha^{-\frac{1}{p}}x)-\big(\alpha^{-\frac{2}{p}}
\lambda(\alpha)-\hat{\lambda}\big)\Big|\hat{u}\\
&+C\alpha^{-\frac{l+1}{p}}\big(1+|x|^{2p+l-1}\big)
\big(\tilde{w}_{\alpha}
+\hat{u}\big)\\
\leq& C\delta_{\alpha}\alpha^{-\frac{2}{p}}e^{-\frac{|x|}{4}} \ \ \text{in}\ \ B_{r_{0}\alpha^{\frac{1}{2p}}}(0)\  \ \text{as}\ \ \alpha\rightarrow\infty,
\end{split}
\end{equation}
where $r_{0}>0$ and $c_{0}>0$ are as in ($M$), and $\delta_{\alpha}>0$ satisfies $\delta_{\alpha}=o(1)$ as $\alpha\rightarrow\infty$.

We now claim that there exists a constant $C>0$, independent of $\alpha>0$, such that
\begin{equation}\label{3.23}
|\eta_{\alpha}(x)|\leq C\delta_{\alpha}\alpha^{-\frac{2}{p}}
\ \ \text{uniformly in}\ \ \mathbb{R}^{N} \ \ \text{as}\ \ \alpha\rightarrow\infty,
\end{equation}
where $\delta_{\alpha}>0$ satisfies $\delta_{\alpha}=o(1)$ as $\alpha\rightarrow\infty$. On the contrary, suppose the above claim (\ref{3.23}) is false. Then up to a subsequence if necessary, there exists a positive constant $C>0$ such that
\begin{equation}\label{3.24}
\frac{\|\eta_{\alpha}(x)\|_{L^{\infty}(\mathbb{R}^{N})}}
{\alpha^{-\frac{2}{p}}}\geq C
\ \ \text{as}\ \ \alpha\rightarrow\infty.
\end{equation}
Denote $\bar{\eta}_{\alpha}:=\frac{\eta_{\alpha}}
{\|\eta_{\alpha}\|_{L^{\infty}(\mathbb{R}^{N})}}$, so that
\begin{equation}\label{3.25}
\|\bar{\eta}_{\alpha}\|_{L^{\infty}(\mathbb{R}^{N})}\equiv 1
\ \ \text{and}\ \ \int_{\mathbb{R}^N}\bar{\eta}_{\alpha}\hat{u}dx\equiv 0
\ \ \text{for any}\ \ \alpha>0,
\end{equation}
due to (\ref{3.58M1}).
By (\ref{3.47})--(\ref{3.13M}), we then derive from (\ref{2.9}) and (\ref{3.24}) that for sufficiently large $\alpha>0$,
\begin{equation}\label{3.53}
\begin{split}
|\bar{\eta}_{\alpha}(x)|=&\frac{
\Big|\tilde{w}_{\alpha}(x)-\big(\int_{\mathbb{R}^{N}}\tilde{w}_{\alpha}
\hat{u}dx\big)\hat{u}(x)\Big|}
{\|\eta_{\alpha}\|_{L^{\infty}(\mathbb{R}^{N})}}\\
\leq&\frac{Ce^{-\frac{|x|}{2}}}{|\alpha^{-\frac{1}{2p}}x|^{4}}
\leq C(r_{0})e^{-\frac{|x|}{4}}
\ \ \text{in}\ \ \mathbb{R}^{N}\backslash B_{r_{0}\alpha^{\frac{1}{2p}}}(0),
\end{split}
\end{equation}
and
\begin{equation}\label{3.54}
\begin{split}
|\nabla\bar{\eta}_{\alpha}|
=&\frac{\Big|\nabla\tilde{w}_{\alpha}-\big(\int_{\mathbb{R}^{N}}\tilde{w}_{\alpha}
\hat{u}dx\big)\nabla\hat{u}\Big|}
{\|\eta_{\alpha}\|_{L^{\infty}(\mathbb{R}^{N})}}\\
\leq& \frac{Ce^{-\frac{|x|}{3}}}{|\alpha^{-\frac{1}{2p}}x|^{4}}
\leq C(r_{0})e^{-\frac{|x|}{5}}
\ \ \text{in}\ \ \mathbb{R}^{N}\backslash B_{r_{0}\alpha^{\frac{1}{2p}}}(0),
\end{split}
\end{equation}
where $r_{0}>0$ is as in ($M$).
We also deduce from (\ref{3.59}) and (\ref{3.24}) that
\begin{equation}\label{3.26}
\begin{split}
|\mathcal{L}_{\alpha}\bar{\eta}_{\alpha}|=&
\frac{|\mathcal{L}_{\alpha}\eta_{\alpha}|}
{\|\eta_{\alpha}\|_{L^{\infty}(\mathbb{R}^{N})}}
\leq\frac{C\delta_{\alpha}\alpha^{-\frac{2}{p}}e^{-\frac{|x|}{4}}}
{\|\eta_{\alpha}\|_{L^{\infty}(\mathbb{R}^{N})}}\\
\leq& C\delta_{\alpha}e^{-\frac{|x|}{4}}
 \ \ \,\text{in}\ \ B_{r_{0}\alpha^{\frac{1}{2p}}}(0)\ \ \text{as}\ \ \alpha\rightarrow\infty,
\end{split}
\end{equation}
where $r_{0}>0$ is as in ($M$).
Applying the comparison principle, one can deduce from (\ref{3.53}) and (\ref{3.26}) that there exist constants $R_{2}>0$ and $C>0$, independent of $\alpha>0$, such that
\begin{equation}\label{3.55}
|\bar{\eta}_{\alpha}(x)|\leq Ce^{-\frac{|x|}{6}} \ \ \,\text{in}\ \ B_{r_{0}\alpha^{\frac{1}{2p}}}(0)\backslash B_{R_{2}}(0)
\ \ \text{as}\ \ \alpha\rightarrow\infty,
\end{equation}
where $r_{0}>0$ is as in ($M$).
Suppose $y_\alpha$ is a global maximum point of $|\bar{\eta}_{\alpha}(x)|$, so that $$|\bar{\eta}_{\alpha}(y_\alpha)|=\max_{x\in \mathbb{R}^{N}}\frac{|\eta_{\alpha}(x)|}
{\|\eta_{\alpha}\|_{L^{\infty}(\mathbb{R}^{N})}}=1.$$
It then follows from (\ref{3.53}) and (\ref{3.55}) that $|y_\alpha|\leq C$ uniformly in $\alpha>0$.

On the other hand,
define
\begin{equation*}
\begin{split}
A_{\alpha}(x):=&p(N+p-2)g(\alpha^{-\frac{1}{p}}x)|x|^{p-2}+
\frac{p^{2}}{\epsilon}g^{2}(\alpha^{-\frac{1}{p}}x)
|x|^{2p-2}\\
&+
\alpha^{-\frac{2}{p}}V(\alpha^{-\frac{1}{p}}x)
-\alpha^{-\frac{2}{p}}\lambda(\alpha),
\ \ x\in B_{r_{0}\alpha^{\frac{1}{2p}}}(0),
\end{split}
\end{equation*}
so that $\mathcal{L}_{\alpha}=-\epsilon\Delta+A_{\alpha}(x)$  in $B_{r_{0}\alpha^{\frac{1}{2p}}}(0)$,
where $r_{0}>0$ is as in ($M$).
We then deduce from (\ref{3.53}) and (\ref{3.26}) that there exists a constant $C>0$, independent of $\alpha>0$, such that
for $r_{0}>0$ given by ($M$),
\begin{equation}\label{3.28}
\begin{split}
\|\bar{\eta}_{\alpha}\|_{H^{1}\big(B_{r_{0}\alpha^{\frac{1}{2p}}}(0)\big)}^{2}
=&\int_{B_{r_{0}\alpha^{\frac{1}{2p}}}(0)}|\nabla \bar{\eta}_{\alpha}|^{2}dx
+\int_{B_{r_{0}\alpha^{\frac{1}{2p}}}(0)}|\bar{\eta}_{\alpha}|^{2}dx\\
\leq&\max\Big\{\frac{1}{\epsilon}, 1\Big\}\Big[\epsilon\int_{B_{r_{0}\alpha^{\frac{1}{2p}}}(0)}|\nabla \bar{\eta}_{\alpha}|^{2}dx+\int_{B_{R_3}(0)}|\bar{\eta}_{\alpha}|^{2}dx\\
&\qquad\qquad\,\quad+\int_{B_{r_{0}\alpha^{\frac{1}{2p}}}(0)\backslash B_{R_3}(0)}A_{\alpha}(x)|\bar{\eta}_{\alpha}|^{2}dx\Big]\\
\leq&\max\Big\{\frac{1}{\epsilon}, 1\Big\}\Big[\int_{B_{r_{0}\alpha^{\frac{1}{2p}}}(0)}(\mathcal{L}_{\alpha}\bar{\eta}_{\alpha})\bar{\eta}_{\alpha}dx
+C(R_{3})\Big]\\
\leq& C\ \ \text{as}\ \ \alpha\rightarrow\infty,
\end{split}
\end{equation}
where $R_{3}>0$ is independent of $\alpha>0$ and chosen such that $A_{\alpha}(x)>1$ if $|x|>R_{3}$.
Furthermore, we obtain from (\ref{3.53}) and (\ref{3.54}) that there exists a constant $C>0$, independent of $\alpha>0$, such that
\begin{equation}\label{3.29}
\|\bar{\eta}_{\alpha}\|_{H^{1}\big(\mathbb{R}^{N}\backslash B_{r_{0}\alpha^{\frac{1}{2p}}}(0)\big)}^{2}\leq C\int_{\mathbb{R}^{N}}e^{-\frac{2|x|}{5}}dx\leq C
\ \ \text{as}\ \ \alpha\rightarrow\infty,
\end{equation}
where $r_{0}>0$ is as in ($M$).

We thus derive from (\ref{3.28}) and (\ref{3.29}) that the sequence $\{\bar{\eta}_{\alpha}\}$ is bounded uniformly in $H^{1}(\mathbb{R}^{N})$ as $\alpha\rightarrow\infty$. Therefore, up to a subsequence if necessary, we may assume that there exists a function $\bar{\eta}_{0}\in H^1(\mathbb{R}^{N})$ such that $$\bar{\eta}_{\alpha}\rightharpoonup\bar{\eta}_{0}\ \ \text{weakly in}\ \ H^{1}(\mathbb{R}^{N})\ \ \text{as}\ \ \alpha\rightarrow\infty,$$
and
$$\bar{\eta}_{\alpha}\rightarrow\bar{\eta}_{0}\ \ \text{strongly in}\ \ L^{2}_{loc}(\mathbb{R}^{N})
\ \ \text{as}\ \ \alpha\rightarrow\infty.$$
It then follows from (\ref{3.25}) and (\ref{3.26}) that $\bar{\eta}_{0}$ satisfies
$$\int_{\mathbb{R}^N}\bar{\eta}_{0}\hat{u}dx=0\ \ \text{and}\ \  \mathcal{L}\bar{\eta}_{0}=0\ \ \text{in}\ \ \mathbb{R}^{N}.$$
We hence derive from (\ref{A2}) and $\int_{\mathbb{R}^N}\hat{u}^{2}dx=1$ that $\bar{\eta}_{0}(x)\equiv0$ in $\mathbb{R}^{N}$, which however contradicts to the fact that up to a subsequence if necessary, $1\equiv|\bar{\eta}_{\alpha}(y_\alpha)|\rightarrow|\bar{\eta}_{0}(\bar{y}_0)|$ as $\alpha\rightarrow\infty$ for some $\bar{y}_0\in \mathbb{R}^{N}$. Thus, the claim (\ref{3.23}) holds true, which implies from (\ref{3.58}) that
\begin{equation}\label{3.201M}
\tilde{w}_{\alpha}(x)
=\Big(\int_{\mathbb{R}^{N}}\tilde{w}_{\alpha}\hat{u}dx\Big)
\hat{u}(x)+o(\alpha^{-\frac{2}{p}})\ \ \text{in}\ \ \mathbb{R}^{N}\ \ \text{as}\ \ \alpha\rightarrow\infty.
\end{equation}
The proof of Lemma \ref{lem3.2} is therefore complete.
\qed

We emphasize that the construction of $\eta_{\alpha}$ in (\ref{3.58}) is crucial for proving Lemma \ref{lem3.2}, which needs the orthogonality (\ref{3.58M1}) between $\eta_{\alpha}(x)$ and $\hat{u}(x)$. For this reason, one cannot replace the coefficient term $\int_{\mathbb{R}^{N}}\tilde{w}_{\alpha}\hat{u}dx$ of (\ref{3.58}) by the constant $1$, in spite of the fact that $\int_{\mathbb{R}^{N}}\tilde{w}_{\alpha}\hat{u}dx$ approaches to $1$ as $\alpha\to\infty$. Based on Lemma \ref{lem3.2}, we next utilize fully the $L^{2}$--constraint conditions of the principal eigenfunction $u_{\alpha}$ to improve the estimate of Lemma \ref{lem3.2}.

\begin{lem}\label{lem3.21}
Under the assumptions of Lemma \ref{lem3.1}, the principal eigenfunction $u_{\alpha}$ of (\ref{1.1}) satisfies
\begin{equation}\label{3.4}
\tilde{w}_{\alpha}(x):=\alpha^{-\frac{N}{2p}}u_{\alpha}
(\alpha^{-\frac{1}{p}}x)
=\hat{u}(x)+o(\alpha^{-\frac{2}{p}})
\ \ \text{in}\ \ \mathbb{R}^{N} \  \text{as}\ \ \alpha\rightarrow\infty,
\end{equation}
where $u_\alpha(x)\equiv0$ in $\mathbb{R}^{N}\backslash\Omega$, and $\hat{u}(x)>0$ is the unique principal eigenfunction of  $\hat{\lambda}$ defined by (\ref{2.3}).
\end{lem}

\noindent\textbf{Proof.}
In view of (\ref{3.47}) and (\ref{3.201}), we define
\begin{equation}\label{3.40}
\tilde{w}_{\alpha}(x):=\hat{u}(x)+\kappa_{\alpha}(x)
+o(\alpha^{-\frac{2}{p}})
\ \ \text{in}\ \ \mathbb{R}^N\   \text{as}\ \ \alpha\rightarrow\infty,
\end{equation}
so that $\kappa_{\alpha}(x)$ satisfies
\begin{equation}\label{3.39M}
\int_{\mathbb{R}^{N}}\kappa_{\alpha}(x)\hat{u}(x)dx\rightarrow0,
\ \,\
\kappa_{\alpha}(x)\rightarrow 0 \ \ \text{uniformly in}\ \ L^{\infty}(\mathbb{R}^N) \ \ \text{as}\ \ \alpha\rightarrow\infty.
\end{equation}
Substituting (\ref{3.40}) into (\ref{3.201}) yields that
\begin{equation}\label{3.41}
\kappa_{\alpha}(x)=\Big(\int_{\mathbb{R}^N}\kappa_{\alpha}\hat{u}dx\Big)\hat{u}(x)
+o(\alpha^{-\frac{2}{p}})
\ \ \text{in}\ \ \mathbb{R}^N\  \text{as}\ \ \alpha\rightarrow\infty.
\end{equation}

In addition, multiplying (\ref{3.201}) by $\tilde{w}_{\alpha}(x)$ and integrating over $\mathbb{R}^N$, we derive from (\ref{3.40}) that
\begin{equation}\label{3.31}
\begin{split}
1=\int_{\mathbb{R}^{N}}\tilde{w}^{2}_{\alpha}dx
=&\Big(\int_{\mathbb{R}^{N}}\tilde{w}_{\alpha}\hat{u}dx\Big)^{2}
+o(\alpha^{-\frac{2}{p}})\int_{\mathbb{R}^{N}}
\tilde{w}_{\alpha}dx\\
=&\Big(1+\int_{\mathbb{R}^{N}}\kappa_{\alpha}\hat{u}dx
+o(\alpha^{-\frac{2}{p}})
\Big)^{2}
+o(\alpha^{-\frac{2}{p}})\ \ \text{as}\ \ \alpha\rightarrow\infty,
\end{split}
\end{equation}
due to the fact that
$\int_{\mathbb{R}^{N}}\hat{u}^{2}dx=1.$
It then follows from (\ref{3.31}) that
\begin{equation*}
\Big(\int_{\mathbb{R}^{N}}\kappa_{\alpha}\hat{u}dx\Big)^{2}
+2\int_{\mathbb{R}^{N}}\kappa_{\alpha}\hat{u}dx
+o(\alpha^{-\frac{2}{p}})=0
\ \ \text{as}\ \ \alpha\rightarrow\infty,
\end{equation*}
which implies from (\ref{3.39M}) that
\begin{equation}\label{3.44}
\int_{\mathbb{R}^N}\kappa_{\alpha}\hat{u}dx=o(\alpha^{-\frac{2}{p}})\ \ \text{as}\ \ \alpha\rightarrow\infty.
\end{equation}
We thus obtain from (\ref{3.41}) and (\ref{3.44}) that
\begin{equation}\label{0.3}
\kappa_{\alpha}(x)=o(\alpha^{-\frac{2}{p}})
\ \ \text{in}\ \ \mathbb{R}^N\ \ \text{as}\ \ \alpha\rightarrow\infty.
\end{equation}
Finally, we conclude from (\ref{3.40}) and (\ref{0.3}) that (\ref{3.4}) holds true,
which hence completes the proof of Lemma \ref{lem3.21}.
\qed

\subsection{Proof of Theorem \ref{cor 1.2}}
In this subsection, we shall finish the proof of  Theorem \ref{cor 1.2}.
Towards this purpose, we also need  the following existence and uniqueness of solutions.

\begin{lem}\label{lem:A.1} For any fixed $\epsilon>0$ and $p\in\{2\}\cup(3,\infty)$,
suppose that $(\hat{\lambda}, \hat{u})$ is the principal eigenpair of (\ref{2.3}), and $f(x)\in L^{2}(\R^{N})\cap C^{\gamma_{2}}(\mathbb{R}^{N})~(0<\gamma_{2}<1)$ satisfies $\int_{\mathbb{R}^{N}}f(x)\hat{u}dx=0$.
Then  for any $c_{0}> 0$,
\begin{equation}\label{A3}
\left\{\begin{split}
&-\epsilon\Delta u+\Big[\frac{p^{2}c_{0}^{2}}{\epsilon}
|x|^{2p-2}+pc_{0}(N+p-2)|x|^{p-2}-\hat{\lambda}\Big]u=f,\\[2mm]
&\int_{\mathbb{R}^{N}}u\hat{u}dx=0
\end{split}\right.
\end{equation}
has a unique solution $\hat{\psi}\in C^{2}(\mathbb{R}^{N})$.
\end{lem}

\noindent\textbf{Proof.}
Recall from (\ref{4}) that
$$\mathcal{L}:=-\epsilon\Delta+\Big[\frac{p^{2}c_{0}^{2}}{\epsilon}
|x|^{2p-2}+pc_{0}(N+p-2)|x|^{p-2}-\hat{\lambda}\Big]\ \ \text{in}\ \ \mathbb{R}^{N}.$$
For convenience, we also denote
$$\mathcal{L}_{0}:=-\epsilon\Delta+\Big[\frac{p^{2}c_{0}^{2}}
{\epsilon}|x|^{2p-2}+pc_{0}(N+p-2)|x|^{p-2}+1\Big]\ \ \text{in}\ \ \mathbb{R}^{N}.$$
It then follows from (\ref{A2}) that $\mathcal{L}_{0}>\mathcal{L}\geq0$, and thus $ker\mathcal{L}_{0}=span\{0\}$. It hence yields from \cite[Theorem 4.1]{O} that $\mathcal{L}_{0}^{-1}$ exists and is a compact operator on $L^2(\mathbb{R}^{N})$. One can also note from (\ref{A3}) that $\mathcal{L}u=f$ admits solutions, if and only if $\mathcal{L}_{0}u-(1+\hat{\lambda})u=f$ admits solutions. Moreover, we deduce from (\ref{A2}) that
$$ker\Big(I-\Big(\big(1+\hat{\lambda}\big)\mathcal{L}_{0}^{-1}
\Big)^*\Big)=span\big\{\hat{u}\big\}.$$
Since
$$\langle\mathcal{L}_{0}^{-1}f,\hat{u}\rangle=
\Big\langle\mathcal{L}_{0}^{-1}f,\frac{\mathcal{L}_{0}\hat{u}}
{1+\hat{\lambda}}\Big\rangle=\Big\langle f,\frac{\hat{u}}
{1+\hat{\lambda}}\Big\rangle=\frac{1}{1+\hat{\lambda}}
\int_{\mathbb{R}^{N}}f(x)\hat{u}dx=0,$$
we now conclude from the Riesz-Schauder theory of compact linear operators that the first equation
of (\ref{A3}) admits a solution $\hat{\psi}_{0}\in C^{2}(\mathbb{R}^{N})$,
where  the standard elliptic regularity theory is also employed. Furthermore, one can check that the function
$$\hat{\psi}(x):=\hat{\psi}_{0}(x)-\Big(\int_{\mathbb{R}^{N}}\hat{\psi}_{0}\hat{u}dx\Big)
\hat{u}(x)\in C^{2}(\mathbb{R}^{N})$$
solves (\ref{A3}), where $\int_{\mathbb{R}^{N}}\hat{u}^{2}dx=1$ is used. This proves the existence of solutions for (\ref{A3}).

We next prove the uniqueness of solutions for (\ref{A3}). On the contrary, suppose $\hat{\psi}$ and $\bar{\psi}$ are two different solutions of (\ref{A3}). We then obtain from  (\ref{A3}) that $\mathcal{L}(\hat{\psi}-\bar{\psi})=0$.
We thus derive from (\ref{A2}) that $\hat{\psi}\equiv\bar{\psi}+c_{2}\hat{u}$ holds in $\mathbb{R}^{N}$ for some constant $c_{2}\in\R$. The second equation of (\ref{A3}) further yields that $c_{2}=0$, which  hence proves that $\hat{\psi}\equiv\bar{\psi}$ in $\mathbb{R}^{N}$.
This completes the proof of Lemma \ref{lem:A.1}.
\qed

One can check from Lemma \ref{lem:A.1} that (\ref{1.121}) and (\ref{1.191}) admit a unique solution $\hat{\psi}_{i}$ for $i=1, 2$, and (\ref{1.12}) and (\ref{1.19}) admit a unique solution $\hat{\psi}_{i}$ for $i=3, 4$.
Applying the previous several lemmas, we are now ready to  complete the proof of Theorem \ref{cor 1.2}.
\vskip 0.05truein

\noindent\textbf{Proof of Theorem \ref{cor 1.2}.} Following (\ref{3.2}), we define
\begin{equation}\label{3.57}
\alpha^{-\frac{2}{p}}\lambda(\alpha):=\hat{\lambda}+\alpha^{-\frac{2}{p}}V(0)+\theta_{\alpha}+o(\theta_{\alpha})
\ \ \text{as}\ \ \alpha\rightarrow\infty,
\end{equation}
where $\theta_{\alpha}$ satisfies $\lim_{\alpha\rightarrow\infty}\frac{\theta_{\alpha}}{\alpha^{-\frac{2}{p}}}=0$.
As before, we also set
\begin{equation*}
\tilde{v}_\alpha(x):=\tilde{w}_\alpha(x)-\hat{u}(x),\ \ x\in\mathbb{R}^{N},
\end{equation*}
where $\tilde{w}_\alpha(x)\geq0$ is defined by (\ref{3.42}), and $\hat{u}(x)>0$ is the unique principal eigenfunction of  $\hat{\lambda}$ defined by (\ref{2.3}). We then deduce that $\tilde{v}_\alpha$ satisfies (\ref{3.12}).
It follows from (\ref{3.12}), (\ref{3.15}), (\ref{3.14}), (\ref{3.4}) and (\ref{3.57}) that for $r_{0}>0$ and $l>3$ given by ($M$),
\begin{align}
o(e^{-r_{0}\alpha^{\frac{1}{2p}}})
=&\int_{B_{r_{0}\alpha^{\frac{1}{2p}}}(0)}(\mathcal{L}\tilde{v}_{\alpha})\hat{u}dx\nonumber\\
=&\int_{B_{r_{0}\alpha^{\frac{1}{2p}}}(0)}B_{\alpha}(x)\tilde{w}_{\alpha}\hat{u}dx\label{3.56}\\
=&\int_{B_{r_{0}\alpha^{\frac{1}{2p}}}(0)}
\Big[-\alpha^{-\frac{2}{p}}V(\alpha^{-\frac{1}{p}}x)
+\big(\alpha^{-\frac{2}{p}}\lambda(\alpha)
-\hat{\lambda}\big)\Big]\tilde{w}_{\alpha}\hat{u}dx+O(\alpha^{-\frac{l+1}{p}})\nonumber\\
=&-\alpha^{-\frac{3}{p}}\int_{\mathbb{R}^{N}}\big [x\cdot\nabla V(0)\big]\hat{u}^{2}dx
+o(\alpha^{-\frac{3}{p}})+\theta_{\alpha}+o(\theta_{\alpha})\ \ \text{as}\ \ \alpha\rightarrow\infty,\nonumber
\end{align}
where the term $B_{\alpha}(x)$ is as in (\ref{3.12}).
We then obtain from (\ref{3.56}) that
\begin{equation*}
\theta_{\alpha}=\alpha^{-\frac{3}{p}}\int_{\mathbb{R}^{N}}\big(x\cdot\nabla V(0)\big)\hat{u}^{2}dx
+o(\alpha^{-\frac{3}{p}})=o(\alpha^{-\frac{3}{p}})\ \ \text{as}\ \ \alpha\rightarrow\infty,
\end{equation*}
due to the radial symmetry of $\hat{u}$.
Combining this with (\ref{3.57}), it yields that
\begin{equation}\label{3.61}
\alpha^{-\frac{2}{p}}\lambda(\alpha)=\hat{\lambda}+\alpha^{-\frac{2}{p}}V(0)+o(\alpha^{-\frac{3}{p}})
\ \ \text{as}\ \ \alpha\rightarrow\infty.
\end{equation}

We now calculate from (\ref{0.8})--(\ref{0.10}), (\ref{0.5}), (\ref{0.6}) and (\ref{3.61}) that for $r_{0}>0$ and $l>3$ given by ($M$),
\begin{equation}\label{3.7}
\begin{split}
&\mathcal{L}_{\alpha}\big(\tilde{w}_{\alpha}
-\hat{u}\big)\\
=&-\Big[|x|^{p}\Delta g(\alpha^{-\frac{1}{p}}x)+2p|x|^{p-2}\big(x\cdot\nabla g(\alpha^{-\frac{1}{p}}x)\big)\\
&+\frac{1}{\epsilon}|x|^{2p}|\nabla g(\alpha^{-\frac{1}{p}}x)|^{2}+\frac{p}{\epsilon}|x|^{2p-2}
\Big(x\cdot\nabla g^{2}(\alpha^{-\frac{1}{p}}x)\Big)\Big]\tilde{w}_{\alpha}
\\
&-\Big[
\alpha^{-\frac{2}{p}}V(\alpha^{-\frac{1}{p}}x)
+\frac{p^{2}}{\epsilon}|x|^{2p-2}\big(g^{2}(\alpha^{-\frac{1}{p}}x)-c^{2}_0\big)\\
&+p(N+p-2)|x|^{p-2}\big(g(\alpha^{-\frac{1}{p}}x)-c_0\big)
-\big(\alpha^{-\frac{2}{p}}
\lambda(\alpha)-\hat{\lambda}\big)
\Big]\hat{u}\\
=&-\Big[\alpha^{-\frac{2}{p}}V(\alpha^{-\frac{1}{p}}x)-\big(\alpha^{-\frac{2}{p}}
\lambda(\alpha)-\hat{\lambda}\big)\Big]\hat{u}
+o(\alpha^{-\frac{3}{p}})\\
=& -\alpha^{-\frac{3}{p}}\big(x\cdot\nabla V(0)\big)\hat{u}+o(\alpha^{-\frac{3}{p}}) \ \ \text{in}\ \ B_{r_{0}\alpha^{\frac{1}{2p}}}(0)\  \ \,\text{as}\ \ \alpha\rightarrow\infty.
\end{split}
\end{equation}
Following (\ref{3.7}), we define
\begin{equation}\label{3.37M}
\tilde{\eta}_{\alpha}(x):=\tilde{w}_{\alpha}(x)-\Big(\int_{\mathbb{R}^{N}}\tilde{w}_{\alpha}\hat{u}dx\Big)
\hat{u}(x)-\alpha^{-\frac{3}{p}}\hat{\psi}_{1}(x)\ \ \text{in}\ \ \mathbb{R}^{N},
\end{equation}
where $\hat{\psi}_{1}$ is given uniquely by
\begin{equation}\label{2M}
\left\{\begin{split}
&-\epsilon\Delta \hat{\psi}_{1}+\Big[\frac{p^{2}c_{0}^{2}}{\epsilon}|x|^{2p-2}
+pc_{0}(N+p-2)|x|^{p-2}-\hat{\lambda}\Big]\hat{\psi}_{1}\\
&=-\big(x\cdot\nabla V(0)\big)\hat{u}\  \ \text{in}\ \ \mathbb{R}^{N},\\
&\int_{\mathbb{R}^{N}}\hat{\psi}_{1}\hat{u}dx=0,\qquad\qquad\qquad
\end{split}\right.
\end{equation}
due to  Lemma \ref{lem:A.1}.
Applying the comparison principle to (\ref{2M}), we can derive from (\ref{2.9}) that there exists a constant $C>0$ such that
$$ |\hat{\psi}_{1}(x)|\leq Ce^{-\frac{|x|}{3}} \ \ \text{in}\ \ \mathbb{R}^{N}.$$
Furthermore, we obtain from (\ref{3.58M}) and (\ref{2M}) that
\begin{equation}\label{0.12}
\begin{split}
\mathcal{L}_{\alpha}\hat{\psi}_{1}=&\Big\{
\frac{p^{2}}{\epsilon}|x|^{2p-2}\big[g^{2}(\alpha^{-\frac{1}{p}}x)-c^{2}_0\big]
+p(N+p-2)|x|^{p-2}\big[g(\alpha^{-\frac{1}{p}}x)-c_0\big]\\
&+\alpha^{-\frac{2}{p}}V(\alpha^{-\frac{1}{p}}x)
-\big[\alpha^{-\frac{2}{p}}
\lambda(\alpha)-\hat{\lambda}\big]
\Big\}\hat{\psi}_{1}-\big(x\cdot\nabla V(0)\big)\hat{u}\ \ \text{in}\ \ \Omega'_{\alpha},
\end{split}
\end{equation}
where the domain $\Omega'_{\alpha}$ is as in (\ref{3.42}).

It hence follows from (\ref{0.8})--(\ref{0.10}), (\ref{0.5}), (\ref{0.6}), (\ref{3.4}), (\ref{3.61}) and (\ref{0.12}) that there exists a constant $C>0$, independent of $\alpha>0$, such that
for $r_{0}>0$ and $l>3$ given by ($M$),
\begin{equation}\label{3.5}
\begin{split}
|\mathcal{L}_{\alpha}\tilde{\eta}_{\alpha}|
=&\Big|\mathcal{L}_{\alpha}\tilde{w}_{\alpha}
-\Big(\int_{\mathbb{R}^{N}}\tilde{w}_{\alpha}\hat{u}dx\Big)
\mathcal{L}_{\alpha}\hat{u}-\alpha^{-\frac{3}{p}}\mathcal{L}_{\alpha}\hat{\psi}_{1}\Big|\\
=&\Big|-\Big[|x|^{p}\Delta g(\alpha^{-\frac{1}{p}}x)+2p|x|^{p-2}\big(x\cdot\nabla g(\alpha^{-\frac{1}{p}}x)\big)+\frac{1}{\epsilon}|x|^{2p}|\nabla g(\alpha^{-\frac{1}{p}}x)|^{2}\\
&+\frac{p}{\epsilon}|x|^{2p-2}
\Big(x\cdot\nabla g^{2}(\alpha^{-\frac{1}{p}}x)\Big)\Big]\tilde{w}_{\alpha}
-\Big[\Big(\int_{\mathbb{R}^{N}}\tilde{w}_{\alpha}\hat{u}dx\Big)\hat{u}
+\alpha^{-\frac{3}{p}}\hat{\psi}_{1}\Big]\\
&
\cdot\Big[
\frac{p^{2}}{\epsilon}|x|^{2p-2}\big(g^{2}(\alpha^{-\frac{1}{p}}x)-c^{2}_0\big)
+p(N+p-2)|x|^{p-2}\big(g(\alpha^{-\frac{1}{p}}x)-c_0\big)\\
&+\alpha^{-\frac{2}{p}}V(\alpha^{-\frac{1}{p}}x)
-\big(\alpha^{-\frac{2}{p}}
\lambda(\alpha)-\hat{\lambda}\big)
\Big]+\alpha^{-\frac{3}{p}}\big(x\cdot\nabla V(0)\big)\hat{u}\Big|\\
\leq& \Big|-\Big(\int_{\mathbb{R}^{N}}\tilde{w}_{\alpha}\hat{u}dx\Big)
\Big[\alpha^{-\frac{2}{p}}V(\alpha^{-\frac{1}{p}}x)-\big(\alpha^{-\frac{2}{p}}
\lambda(\alpha)-\hat{\lambda}\big)\Big]\hat{u}\\
&-\alpha^{-\frac{3}{p}}
\Big[\alpha^{-\frac{2}{p}}V(\alpha^{-\frac{1}{p}}x)-\big(\alpha^{-\frac{2}{p}}
\lambda(\alpha)-\hat{\lambda}\big)\Big]\hat{\psi}_{1}
+\alpha^{-\frac{3}{p}}\big(x\cdot\nabla V(0)\big)\hat{u}\Big|\\
&+C\big(1+|x|^{2p+l-1}\big)
\Big(\alpha^{-\frac{l+1}{p}}\tilde{w}_{\alpha}
+\alpha^{-\frac{l+1}{p}}\hat{u}+\alpha^{-\frac{l+4}{p}}|\hat{\psi}_{1}|\Big)\\
=&\Big|-\alpha^{-\frac{3}{p}}\big(1+o(\alpha^{-\frac{2}{p}})\big)\big(x\cdot\nabla V(0)+o(1)\big)\hat{u}
+o(\alpha^{-\frac{3}{p}})\hat{\psi}_{1}\\
&+\alpha^{-\frac{3}{p}}\big(x\cdot\nabla V(0)\big)\hat{u}\Big|
+o(\alpha^{-\frac{3}{p}})\big(1+|x|^{2p+l-1}\big)\big(\tilde{w}_{\alpha}
+\hat{u}+|\hat{\psi}_{1}|\big)\\
\leq& C\delta'_{\alpha}\alpha^{-\frac{3}{p}}e^{-\frac{|x|}{4}} \ \ \text{in}\ \ B_{r_{0}\alpha^{\frac{1}{2p}}}(0)\  \ \text{as}\ \ \alpha\rightarrow\infty,
\end{split}
\end{equation}
where $\delta'_{\alpha}>0$ satisfies $\delta'_{\alpha}=o(1)$ as $\alpha\rightarrow\infty$.
Using the similar argument of (\ref{3.23})--(\ref{3.201M}), one can further get from (\ref{3.5}) that
$$|\tilde{\eta}_{\alpha}|\leq C\delta'_{\alpha}\alpha^{-\frac{3}{p}} \ \ \text{uniformly in}\ \ \mathbb{R}^{N}\  \ \text{as}\ \ \alpha\rightarrow\infty,$$
which then gives from (\ref{3.37M}) that
$$\tilde{w}_{\alpha}(x)=\Big(\int_{\mathbb{R}^{N}}\tilde{w}_{\alpha}\hat{u}dx\Big)
\hat{u}(x)+\alpha^{-\frac{3}{p}}\hat{\psi}_{1}+o(\alpha^{-\frac{3}{p}})\ \ \text{in}\ \ \mathbb{R}^{N}
\ \ \text{as}\ \ \alpha\rightarrow\infty.$$
Similar to (\ref{3.40})--(\ref{0.3}), we thus derive from above that
\begin{equation}\label{3.6}
\tilde{w}_{\alpha}(x)=\hat{u}(x)+\alpha^{-\frac{3}{p}}\hat{\psi}_{1}+o(\alpha^{-\frac{3}{p}})
\ \ \text{in}\ \ \mathbb{R}^{N}
\ \ \text{as}\ \ \alpha\rightarrow\infty,
\end{equation}
where $\hat{\psi}_{1}$ satisfies (\ref{0.12}).

In view (\ref{3.61}), we next define
\begin{equation}\label{3.57M}
\alpha^{-\frac{2}{p}}\lambda(\alpha):=\hat{\lambda}+\alpha^{-\frac{2}{p}}V(0)+\tilde{\theta}_{\alpha}
+o(\tilde{\theta}_{\alpha})
\ \ \text{as}\ \ \alpha\rightarrow\infty,
\end{equation}
where $\tilde{\theta}_{\alpha}$ satisfies $\lim_{\alpha\rightarrow\infty}\frac{\tilde{\theta}_{\alpha}}{\alpha^{-\frac{3}{p}}}=0$. Substituting (\ref{3.6}) and (\ref{3.57M}) into (\ref{3.56}), we then obtain that
$$\tilde{\theta}_{\alpha}=\alpha^{-\frac{4}{p}}\sum_{|\sigma|=2}
\frac{D^{\sigma}V(0)}{\sigma!}\int_{\mathbb{R}^{N}}
x^{\sigma}\hat{u}^{2}dx
+o(\alpha^{-\frac{4}{p}})\ \ \text{as}\ \ \alpha\rightarrow\infty,$$
which gives (\ref{3.341}) in view of  (\ref{3.57M}). Finally, we define
$$\hat{\eta}_{\alpha}(x):=\tilde{w}_{\alpha}(x)-\Big(\int_{\mathbb{R}^{N}}\tilde{w}_{\alpha}\hat{u}dx\Big)
\hat{u}(x)-\alpha^{-\frac{3}{p}}\hat{\psi}_{1}(x)-\alpha^{-\frac{4}{p}}\hat{\psi}_{2}(x)\ \ \text{in}\ \ \mathbb{R}^{N},$$
where $\hat{\psi}_{1}$ and $\hat{\psi}_{2}$ are given uniquely by (\ref{1.121}) and (\ref{1.191}).
Repeating the process of (\ref{3.5})--(\ref{3.6}), one can check that (\ref{3.361}) holds true. This therefore completes the proof of Theorem \ref{cor 1.2}.
\qed

\section{Proof of Theorem \ref{thm:1.3}}

This section is devoted to the proof of Theorem \ref{thm:1.3} on the refined limiting profiles of the principal eigenpair for (\ref{1.1}) as $\alpha\rightarrow\infty$, where the regularity of $V(x)$ near the origin is weaker than that of Theorem \ref{cor 1.2}.

Under the assumptions of Theorem \ref{thm:1.3}, we first note that all results of Section 2 are applicable. Similar to Lemmas \ref{lem3.1} and \ref{lem3.21}, we now investigate the following second-order asymptotic expansions of the principal eigenpair $(\lambda(\alpha), u_{\alpha})$ for (\ref{1.1}) as $\alpha\rightarrow\infty$.

\begin{lem}\label{lem:3.01}
Suppose that $0\leq V(x)\in C^{\gamma}(\bar{\Omega})~(0<\gamma<1)$ satisfies  $V(x)=\hat{h}(x)+o(|x|^{\hat{q}})$ as $|x|\rightarrow0$, where $\hat{h}(x)$ satisfies (\ref{1.11A}) for some $\hat{q}>0$. Assume that $m(x)=g(x)|x|^p\geq0$, where $p\in\{2\}\cup(3,\infty)$ and $g(x)$ satisfies ($M$) with $l>\hat{q}+1$.
Then for any fixed $\epsilon>0$, the principal eigenpair $(\lambda(\alpha), u_{\alpha})$ of (\ref{1.1}) satisfies
\begin{equation}\label{1.10}
\lambda(\alpha)=\alpha^{\frac{2}{p}}\hat{\lambda}+\alpha^{-\frac{\hat{q}}{p}}
\int_{\mathbb{R}^{N}}\hat{h}(x)\hat{u}^{2}dx+o(\alpha^{-\frac{\hat{q}}{p}})\ \ \text{as}\ \ \alpha\rightarrow\infty,
\end{equation}
and
\begin{equation}\label{3.20}
\begin{split}
\tilde{w}_{\alpha}(x):=&\alpha^{-\frac{N}{2p}}u_{\alpha}
(\alpha^{-\frac{1}{p}}x)\\
=&\hat{u}(x)+\alpha^{-\frac{\hat{q}+2}{p}}
\hat{\psi}_{3}(x)
+o(\alpha^{-\frac{\hat{q}+2}{p}})
\ \ \text{in}\ \ \mathbb{R}^{N} \  \text{as}\ \ \alpha\rightarrow\infty,
\end{split}
\end{equation}
where $(\hat{\lambda}, \hat{u})$ is the unique principal eigenpair of (\ref{2.3}), $u_\alpha(x)\equiv0$ in $\mathbb{R}^{N}\backslash\Omega$,
and $\hat{\psi}_{3}$ is the unique solution of (\ref{1.12}) and (\ref{1.19}).
\end{lem}

\noindent\textbf{Proof.} Motivated by (\ref{3.7M}), we
define
\begin{equation*}
\tilde{v}_\alpha(x)=\tilde{w}_\alpha(x)-\hat{u}(x)\ \ \text{in} \ \ \mathbb{R}^{N},
\end{equation*}
where $\tilde{w}_\alpha(x)\geq0$ is defined by (\ref{3.20}), so that $\tilde{v}_\alpha$ satisfies $\tilde{v}_\alpha\rightarrow0$ uniformly in $L^{\infty}(\mathbb{R}^{N})$ as $\alpha\rightarrow\infty$ by the argument of (\ref{3.47}). Moreover, $\tilde{v}_\alpha$ satisfies
\begin{equation}\label{3.115M}
\mathcal{L}\tilde{v}_{\alpha}
=B_{\alpha}(x)\tilde{w}_{\alpha}(x)
\ \ \text{in}\ \ \Omega'_{\alpha},
\end{equation}
where the operator $\mathcal{L}$ is as in (\ref{4}), the term $B_{\alpha}(x)$ is as in (\ref{3.12}), and the domain $\Omega'_{\alpha}$ is defined by (\ref{3.42}).  Similar to (\ref{3.56}),
one can get from (\ref{3.115M}) that for sufficiently large $\alp >0,$
\begin{align}
o(e^{-r_{0}\alpha^{\frac{1}{2p}}})
=&\int_{B_{r_{0}\alpha^{\frac{1}{2p}}}(0)}(\mathcal{L}\tilde{v}_{\alpha})\hat{u}dx\nonumber\\
=&\int_{B_{r_{0}\alpha^{\frac{1}{2p}}}(0)}B_{\alpha}(x)\tilde{w}_{\alpha}\hat{u}dx\label{3.115}\\
=&\int_{B_{r_{0}\alpha^{\frac{1}{2p}}}(0)}
\Big[-\alpha^{-\frac{2}{p}}V(\alpha^{-\frac{1}{p}}x)
+\big(\alpha^{-\frac{2}{p}}\lambda(\alpha)
-\hat{\lambda}\big)\Big]\tilde{w}_{\alpha}\hat{u}dx+O(\alpha^{-\frac{l+1}{p}})\nonumber\\
=&-\alpha^{-\frac{\hat{q}+2}{p}}\Big(\int_{\mathbb{R}^{N}}
\hat{h}(x)\hat{u}^{2}dx+o(1)\Big)
+
\big(\alpha^{-\frac{2}{p}}\lambda(\alpha)-\hat{\lambda}\big)\big(1+o(1)\big),\nonumber
\end{align}
where $r_{0}>0$ and $l>\hat{q}+1$ given by ($M$).
This further implies that
\begin{equation*}\label{3.16}
\alpha^{-\frac{2}{p}}\lambda(\alpha)-\hat{\lambda}
=\alpha^{-\frac{\hat{q}+2}{p}}\Big(\int_{\mathbb{R}^{N}}\hat{h}(x)\hat{u}^{2}dx+o(1)\Big)
\ \ \text{as}\ \ \alpha\rightarrow\infty,
\end{equation*}
and thus (\ref{1.10}) holds true.


We next define
\begin{equation}\label{3.21}
\begin{aligned}
v_{\alpha}(x):=
\tilde{w}_{\alpha}(x)-
\Big(\int_{\mathbb{R}^{N}}\tilde{w}_{\alpha}\hat{u}dx\Big)\hat{u}(x)
-\alpha^{-\frac{\hat{q}+2}{p}}\hat{\psi}_{3}(x),
\ \ x\in\mathbb{R}^{N},
\end{aligned}
\end{equation}
where $\hat{\psi}_{3}$ is the unique solution of (\ref{1.12}) and (\ref{1.19}). Note that $\tilde{w}_{\alpha}$ and $\hat{u}$ satisfy (\ref{0.5}) and (\ref{0.6}), respectively.
It then follows from (\ref{1.12}), (\ref{1.19}), (\ref{0.5}) and (\ref{0.6}) that
\begin{equation} \label{3.21M}
\begin{split}
\mathcal{L}_{\alpha}v_{\alpha}=& \mathcal{L}_{\alpha}\tilde{w}_{\alpha}-
\Big(\int_{\mathbb{R}^{N}}\tilde{w}_{\alpha}\hat{u}dx\Big)
\mathcal{L}_{\alpha}\hat{u}
-\alpha^{-\frac{\hat{q}+2}{p}}\mathcal{L}_{\alpha}\hat{\psi}_{3}\\
=&-\Big\{|x|^{p}\Delta g(\alpha^{-\frac{1}{p}}x)+2p|x|^{p-2}\big[x\cdot\nabla g(\alpha^{-\frac{1}{p}}x)\big]+\frac{1}{\epsilon}|x|^{2p}|\nabla g(\alpha^{-\frac{1}{p}}x)|^{2}\\
&+\frac{p}{\epsilon}|x|^{2p-2}
\big[x\cdot\nabla g^{2}(\alpha^{-\frac{1}{p}}x)\big]\Big\}\tilde{w}_{\alpha}
-\Big(\int_{\mathbb{R}^{N}}\tilde{w}_{\alpha}\hat{u}dx\Big)
\Big\{
\alpha^{-\frac{2}{p}}V(\alpha^{-\frac{1}{p}}x)\\
&+p(N+p-2)|x|^{p-2}\big[g(\alpha^{-\frac{1}{p}}x)-c_0\big]
+\frac{p^{2}}{\epsilon}|x|^{2p-2}\big[g^{2}(\alpha^{-\frac{1}{p}}x)-c^{2}_0\big]\\
&
-\big[\alpha^{-\frac{2}{p}}
\lambda(\alpha)-\hat{\lambda}\big]
\Big\}\hat{u}-\alpha^{-\frac{\hat{q}+2}{p}}\Big\{
\frac{p^{2}}{\epsilon}|x|^{2p-2}\big[g^{2}(\alpha^{-\frac{1}{p}}x)-c^{2}_0\big]\hat{\psi}_{3}\\
&+p(N+p-2)|x|^{p-2}\big[g(\alpha^{-\frac{1}{p}}x)-c_0\big]
\hat{\psi}_{3}+\alpha^{-\frac{2}{p}}V(\alpha^{-\frac{1}{p}}x)\hat{\psi}_{3}\\
&
-\big[\alpha^{-\frac{2}{p}}
\lambda(\alpha)-\hat{\lambda}\big]\hat{\psi}_{3}
-\hat{h}(x)\hat{u}+\big(\int_{\mathbb{R}^{N}}\hat{h}(x)\hat{u}^{2}dx\big)
\hat{u}\Big\}
\ \ \mbox{in}\ \ \Omega'_{\alpha},
\end{split}
\end{equation}
where $\Omega'_{\alpha}$ is as in (\ref{3.42}) and the operator $\mathcal{L}_{\alpha}$ is as in (\ref{3.58M}).
Moreover, by the comparison principle, together with Lemma \ref{lem:2.1}, we obtain from (\ref{1.12}) and (\ref{1.19}) that there exists a constant $C>0$ such that
\begin{equation}\label{4.22M}
|\hat{\psi}_{3}(x)|\leq Ce^{-\frac{|x|}{3}}\ \ \text{and}\ \ | \hat{\psi}_{4}(x)|\leq Ce^{-\frac{|x|}{6}}\ \ \text{in}\ \ \mathbb{R}^{N}.
\end{equation}
By applying (\ref{0.8})--(\ref{0.10}), we then calculate from  (\ref{1.10}) and (\ref{3.21M}) that there exists a constant $C>0$, independent of $\alpha>0$, such that for $r_{0}>0$ and $l>\hat{q}+1$ given by ($M$),
\begin{align}
|\mathcal{L}_{\alpha}v_{\alpha}|\leq&
\Big|\Big(\int_{\mathbb{R}^{N}}
\tilde{w}_{\alpha}\hat{u}dx\Big)
\Big[-\alpha^{-\frac{\hat{q}+2}{p}}\hat{h}(x)\hat{u}+o(\alpha^{-\frac{\hat{q}+2}{p}})
+\big(\alpha^{-\frac{2}{p}}\lambda(\alpha)-\hat{\lambda}\big)\hat{u}\Big]\nonumber\\
&-\alpha^{-\frac{\hat{q}+2}{p}}\Big[\Big(\int_{\mathbb{R}^{N}}\hat{h}(x)\hat{u}^{2}dx\Big)
\hat{u}+\alpha^{-\frac{\hat{q}+2}{p}}\hat{h}(x)\hat{\psi}_{3}
+o(\alpha^{-\frac{\hat{q}+2}{p}})\nonumber\\
&\qquad\qquad-\big(\alpha^{-\frac{2}{p}}\lambda(\alpha)-\hat{\lambda}\big)
\hat{\psi}_{3}-\hat{h}(x)\hat{u}\Big]\Big|\label{3.22}\\
&+C\big(1+|x|^{2p+l-1}\big)
\Big(\alpha^{-\frac{l+1}{p}}\tilde{w}_{\alpha}
+\alpha^{-\frac{l+1}{p}}\hat{u}+\alpha^{-\frac{l+\hat{q}+3}{p}}|\hat{\psi}_{3}|\Big)\nonumber\\
\leq&  C\delta_{\alpha}\alpha^{-\frac{\hat{q}+2}{p}}e^{-\frac{|x|}{4}}
\ \ \text{in}\ \ B_{r_{0}\alpha^{\frac{1}{2p}}}(0)
 \ \ \text{as}\ \ \alpha\rightarrow\infty,\nonumber
\end{align}
where $\delta_{\alpha}>0$ satisfies $\delta_{\alpha}=o(1)$ as $\alpha\rightarrow\infty$.
Using the similar argument of (\ref{3.23})--(\ref{3.201M}), one can then derive from  (\ref{3.21}) and (\ref{3.22}) that
\begin{equation}\label{3.30}
\tilde{w}_{\alpha}(x)=\Big(\int_{\mathbb{R}^N}\tilde{w}_{\alpha}\hat{u}dx\Big)\hat{u}(x)
+\alpha^{-\frac{\hat{q}+2}{p}}\hat{\psi}_{3}(x)+o(\alpha^{-\frac{\hat{q}+2}{p}})\  \ \text{in}\ \ \mathbb{R}^N\ \ \text{as}\ \ \alpha\rightarrow\infty.
\end{equation}
Furthermore, similar to (\ref{3.40})--(\ref{0.3}), we thus deduce from (\ref{3.30}) that
\begin{equation}\label{3.8}
\tilde{w}_{\alpha}(x)
=\hat{u}(x)+\alpha^{-\frac{\hat{q}+2}{p}}
\hat{\psi}_{3}(x)
+o(\alpha^{-\frac{\hat{q}+2}{p}})
\ \ \text{in}\ \ \mathbb{R}^{N} \  \text{as}\ \ \alpha\rightarrow\infty,
\end{equation}
and   the proof of Lemma \ref{lem:3.01} is therefore complete.
\qed

Applying Lemma \ref{lem:3.01}, we next complete the proof of Theorem \ref{thm:1.3}.

\vskip 0.05truein
\noindent\textbf{Proof of Theorem \ref{thm:1.3}.}
Following (\ref{1.10}), we first denote
\begin{equation}\label{3.33}
\alpha^{-\frac{2}{p}}\lambda(\alpha)
:=\hat{\lambda}+\alpha^{-\frac{\hat{q}+2}{p}} \int_{\mathbb{R}^{N}}\hat{h}(x)\hat{u}^{2}dx+\beta_{\alpha}
+o(\beta_{\alpha})
\ \ \text{as}\ \ \alpha\rightarrow\infty,
\end{equation}
where $\beta_{\alpha}$ satisfies
$$\lim_{\alpha\rightarrow\infty}\frac{\beta_{\alpha}}
{\alpha^{-\frac{\hat{q}+2}{p}}}=0.$$
Applying Lemma \ref{lem:2.1} and (\ref{4.22M}), substituting (\ref{3.20}) into (\ref{3.115}), we then derive from (\ref{3.33}) that for $r_{0}>0$ and $l>2\hat{q}+3$ given by ($M$),
\begin{align}
&o(e^{-r_{0}\alpha^{\frac{1}{2p}}})\nonumber\\
=&\int_{B_{r_{0}\alpha^{\frac{1}{2p}}}(0)}
\Big[-\alpha^{-\frac{2}{p}}V(\alpha^{-\frac{1}{p}}x)
+\big(\alpha^{-\frac{2}{p}}\lambda(\alpha)
-\hat{\lambda}\big)\Big]\tilde{w}_{\alpha}\hat{u}dx+O(\alpha^{-\frac{l+1}{p}})\nonumber\\
=&\int_{B_{r_{0}\alpha^{\frac{1}{2p}}}(0)}\Big[-\alpha^{-\frac{\hat{q}+2}{p}}\hat{h}(x)
+\alpha^{-\frac{\hat{q}+2}{p}}\int_{\mathbb{R}^{N}}\hat{h}(x)\hat{u}^{2}dx
+\beta_{\alpha}
+o(\beta_{\alpha})\Big]\nonumber\\
&\quad\quad\quad\quad \ \ \cdot\Big(\hat{u}+\alpha^{-\frac{\hat{q}+2}{p}}\hat{\psi}_{3}
+o(\alpha^{-\frac{\hat{q}+2}{p}})\Big)\hat{u}dx+O(\alpha^{-\frac{l+1}{p}})\label{3.35}\\
=&\alpha^{-\frac{\hat{q}+2}{p}}\Big\{-\int_{B_{r_{0}\alpha^{\frac{1}{2p}}}(0)}
\hat{h}(x)\hat{u}^{2}dx+\int_{\mathbb{R}^{N}}\hat{h}(x)\hat{u}^{2}dx
\int_{B_{r_{0}\alpha^{\frac{1}{2p}}}(0)}\hat{u}^{2}dx\Big\}+o(\beta_{\alpha})\nonumber\\
&-\alpha^{-\frac{2\hat{q}+4}{p}}\int_{B_{r_{0}\alpha^{\frac{1}{2p}}}(0)}
\hat{h}(x)\hat{\psi}_{3}\hat{u}dx+\beta_{\alpha}
\int_{B_{r_{0}\alpha^{\frac{1}{2p}}}(0)}\hat{u}^{2}dx
+o(\alpha^{-\frac{2\hat{q}+4}{p}})\nonumber\\
=&-\alpha^{-\frac{2\hat{q}+4}{p}}\int_{\mathbb{R}^{N}}\hat{h}(x)
\hat{\psi}_{3}\hat{u}dx
+\beta_{\alpha}+o(\beta_{\alpha})+o(\alpha^{-\frac{2\hat{q}+4}{p}})
\ \ \text{as}\ \ \alpha\rightarrow\infty,\nonumber
\end{align}
where we have used the identities $\int_{\mathbb{R}^{N}}\hat{u}^{2}dx=1$ and $\int_{\mathbb{R}^{N}}\hat{\psi}_{3}\hat{u}dx=0$, and $\hat{\psi}_{3}$ is given uniquely by (\ref{1.12}) and (\ref{1.19}).
It then yields from (\ref{3.35}) that
$$\beta_{\alpha}=
\alpha^{-\frac{2\hat{q}+4}{p}}\int_{\mathbb{R}^{N}}\hat{h}(x)
\hat{\psi}_{3}\hat{u}dx+o(\alpha^{-\frac{2\hat{q}+4}{p}})\ \ \text{as}\ \ \alpha\rightarrow\infty,$$
which thus proves (\ref{3.34}) in view of (\ref{3.33}).

To prove (\ref{3.36}), we now define
\begin{equation*}\label{3.51}
F_\alpha(x):=\tilde{w}_{\alpha}(x)-\hat{u}(x)-
\alpha^{-\frac{\hat{q}+2}{p}}\hat{\psi}_{3}(x)
\ \ \text{in} \ \ \mathbb{R}^{N},
\end{equation*}
where $\hat{\psi}_{3}$ is given uniquely by (\ref{1.12}) and (\ref{1.19}).
It follows from Lemma \ref{lem:3.01} that
\begin{equation*}
F_\alpha(x)=o\big(\alpha^{-\frac{\hat{q}+2}{p}}\big) \ \,\ \text{uniformly in}\ \, L^{\infty}(\mathbb{R}^N) \ \ \text{as}\ \ \alpha\rightarrow\infty.
\end{equation*}
Applying (\ref{3.34}) and (\ref{3.21M}),
it then follows from (\ref{0.8})--(\ref{0.10}) and (\ref{3.20})  that for $r_{0}>0$ and $l>2\hat{q}+3$ given by ($M$),
\begin{equation}\label{3.52}
\begin{split}
\mathcal{L_{\alpha}}F_\alpha=&\mathcal{L_{\alpha}}\tilde{w}_{\alpha}-\mathcal{L_{\alpha}}\hat{u}-
\alpha^{-\frac{\hat{q}+2}{p}}\mathcal{L_{\alpha}}\hat{\psi}_{3}\\
=&-\Big\{|x|^{p}\Delta g(\alpha^{-\frac{1}{p}}x)+2p|x|^{p-2}\big[x\cdot\nabla g(\alpha^{-\frac{1}{p}}x)\big]+\frac{1}{\epsilon}|x|^{2p}|\nabla g(\alpha^{-\frac{1}{p}}x)|^{2}\\
&+\frac{p}{\epsilon}|x|^{2p-2}
\big[x\cdot\nabla g^{2}(\alpha^{-\frac{1}{p}}x)\big]\Big\}\tilde{w}_{\alpha}
-
\Big\{\frac{p^{2}}{\epsilon}|x|^{2p-2}\big[g^{2}(\alpha^{-\frac{1}{p}}x)-c^{2}_0\big]
\\
&+p(N+p-2)|x|^{p-2}\big[g(\alpha^{-\frac{1}{p}}x)-c_0\big]
-\big[\alpha^{-\frac{2}{p}}
\lambda(\alpha)-\hat{\lambda}\big]
\\
&
+\alpha^{-\frac{2}{p}}V(\alpha^{-\frac{1}{p}}x)
\Big\}\hat{u}-\alpha^{-\frac{\hat{q}+2}{p}}\Big\{
\frac{p^{2}}{\epsilon}|x|^{2p-2}\big[g^{2}(\alpha^{-\frac{1}{p}}x)-c^{2}_0\big]\hat{\psi}_{3}\\
&+p(N+p-2)|x|^{p-2}\big[g(\alpha^{-\frac{1}{p}}x)-c_0\big]
\hat{\psi}_{3}+\alpha^{-\frac{2}{p}}V(\alpha^{-\frac{1}{p}}x)\hat{\psi}_{3}\\
&
-\big[\alpha^{-\frac{2}{p}}
\lambda(\alpha)-\hat{\lambda}\big]\hat{\psi}_{3}
-\hat{h}(x)\hat{u}+\big(\int_{\mathbb{R}^{N}}\hat{h}(x)\hat{u}^{2}dx\big)
\hat{u}\Big\}\\
=&\Big[-\alpha^{-\frac{2}{p}}V(\alpha^{-\frac{1}{p}}x)
+\big(\alpha^{-\frac{2}{p}}\lambda(\alpha)-\hat{\lambda}\big)
\Big]\hat{u}-\alpha^{-\frac{\hat{q}+2}{p}}
\Big[\alpha^{-\frac{2}{p}}V(\alpha^{-\frac{1}{p}}x)\hat{\psi}_{3}\\
&-\big(\alpha^{-\frac{2}{p}}
\lambda(\alpha)-\hat{\lambda}\big)\hat{\psi}_{3}-\hat{h}(x)\hat{u}
+\big(\int_{\mathbb{R}^{N}}\hat{h}(x)\hat{u}^{2}dx\big)
\hat{u}\Big]+o(\alpha^{-\frac{2\hat{q}+4}{p}})\\
=&\alpha^{-\frac{2\hat{q}+4}{p}}\Big[-\hat{h}(x)\hat{\psi}_{3}+
\Big(\int_{\mathbb{R}^{N}}\hat{h}(x)\hat{\psi}_{3}\hat{u}dx\Big)\hat{u}
+\Big(\int_{\mathbb{R}^{N}}\hat{h}(x)\hat{u}^{2}dx\Big)\hat{\psi}_{3}\Big]\\
&
+o(\alpha^{-\frac{2\hat{q}+4}{p}})\ \ \text{in}\ \ B_{r_{0}\alpha^{\frac{1}{2p}}}(0)\,
\ \ \text{as}\ \ \alpha\rightarrow\infty,
\end{split}
\end{equation}
where the operator $\mathcal{L_{\alpha}}$ is defined by (\ref{3.58M}).
In view of (\ref{3.52}), we finally define
\begin{equation}\label{3.38}
M_\alpha(x):=\tilde{w}_{\alpha}(x)-
\Big(\int_{\mathbb{R}^{N}}\tilde{w}_{\alpha}\hat{u}dx\Big)\hat{u}(x)-
\alpha^{-\frac{\hat{q}+2}{p}}\hat{\psi}_{3}(x)
-\alpha^{-\frac{2\hat{q}+4}{p}}\hat{\psi}_{4}(x)
\ \ \text{in} \ \ \mathbb{R}^{N},
\end{equation}
where $\hat{\psi}_{3}(x)$ and $\hat{\psi}_{4}(x)$ are given uniquely by (\ref{1.12}) and (\ref{1.19}).
Repeating the process of (\ref{3.22})--(\ref{3.8}), one can further prove that (\ref{3.36}) holds true. This therefore completes the proof of Theorem \ref{thm:1.3}.\qed

\vskip 0.3truein	
\noindent {\bf Acknowledgements}
\vskip 0.1truein
\noindent The authors thank Professor Rui Peng very much for his fruitful discussions on the present paper.

\vskip 0.3truein	
\noindent {\bf Funding}
\vskip 0.1truein
\noindent Yujin Guo is partially supported by National Key R $\&$ D Program of China (Grant 2023YFA1010001), and NSF of China (Grants 12225106 and 12371113). Yuan Lou is partially supported by NSF of China (Grants 12250710674 and 12261160366).

\vskip 0.3truein

\noindent {\bf Data availability}
\vskip 0.1truein
\noindent No data was used for the research described in the article.

\vskip 0.3truein
\noindent {\bf Disclosure statement}
\vskip 0.1truein
\noindent The authors report that there are no competing interests to declare.

\end{document}